\documentclass[10 pt, journal]{IEEEtran}  
\IEEEoverridecommandlockouts                             

\usepackage{xspace,amssymb,epsfig,subfigure,syntonly}
\usepackage{epsfig,amsmath,color}
\usepackage{xspace,syntonly}
\usepackage{url}
\usepackage{algorithm,algorithmic}
\usepackage{mathtools}
\usepackage{cite}

\newcommand{\utwi}[1]{\mbox{\boldmath $#1$}}

\newcommand{\trace}{{\textrm{Tr}}}

\newcommand{\cD}{{\cal D}}

\newcommand{\cN}{{\cal N}}
\newcommand{\cP}{{\cal P}}

\newcommand{\cE}{{\cal E}}

\newcommand{\cU}{{\cal U}}

\newcommand{\cV}{{\cal V}}

\newcommand{\cY}{{\cal Y}}

\newcommand{\bd}{{\bf d}}
\newcommand{\be}{{\bf e}}
\newcommand{\bbf}{{\bf f}}
\newcommand{\bg}{{\bf g}}
\newcommand{\bh}{{\bf h}}

\newcommand{\bp}{{\bf p}}

\newcommand{\br}{{\bf r}}

\newcommand{\bx}{{\bf x}}
\newcommand{\bu}{{\bf u}}

\newcommand{\bv}{{\bf v}}
\newcommand{\bi}{{\bf i}}

\newcommand{\by}{{\bf y}}

\newcommand{\bL}{{\bf L}}

\newcommand{\bI}{{\bf I}}
\newcommand{\bX}{{\bf X}}
\newcommand{\bZ}{{\bf Z}}

\newcommand{\bY}{{\bf Y}}
\newcommand{\bV}{{\bf V}}

\newcommand{\blambda}{{\utwi{\lambda}}}

\newcommand{\bPsi}{{\utwi{\Psi}}}

\newcommand{\bPhi}{{\utwi{\Phi}}}

\newcommand{\bUpsilon}{{\utwi{\Upsilon}}}

\newcommand{\sfH}{\textsf{H}}
\newcommand{\sfT}{\textsf{T}}

\newcommand{\imj}{\mathrm{j}}

\usepackage{epstopdf}

\begin{document}

\newtheorem{definition}{Definition}
\newtheorem{remark}{Remark}
\newtheorem{proposition}{Proposition}
\newtheorem{lemma}{Lemma}
\newtheorem{assumption}{Assumption}
\newtheorem{theorem}{Theorem}

\markboth{IEEE TRANSACTIONS ON POWER SYSTEMS}%
{ \MakeLowercase{\textit{et al.}}: }


\title{Photovoltaic Inverter Controllers Seeking \\ AC Optimal Power Flow Solutions}

\author{Emiliano Dall'Anese, \emph{Member, IEEE}, Sairaj V. Dhople, \emph{Member, IEEE}, and Georgios B. Giannakis, \emph{Fellow, IEEE}
\thanks{Submitted: December 31, 2014. Revised April 24, 2015 and June 30, 2015. Accepted July 2, 2015.}
\thanks{The work of E. Dall'Anese was supported in part by the Laboratory Directed Research and Development Program at the National Renewable Energy Laboratory, and the Department of Energy Office of Electricity (OE) agreement 28676; S. V. Dhople and G. B. Giannakis were supported in part by the National Science Foundation through grants CCF 1423316, CyberSEES 1442686, and CAREER award ECCS-1453921; and by the Institute of Renewable Energy and the Environment, University of Minnesota under grant RL-0010-13.}
\thanks{E. Dall'Anese is with the National Renewable Energy Laboratory, Golden, CO. S. V. Dhople and G. B. Giannakis are with the Department of Electrical and Computer Engineering, and the Digital Technology Center, University of Minnesota, Minneapolis, MN 55455, USA.   E-mail: {emiliano.dallanese@nrel.gov}, {\{sdhople, georgios\}@umn.edu}.%
}
}
\maketitle

\begin{abstract}

This paper considers future distribution networks featuring inverter-interfaced photovoltaic (PV) systems, and addresses the synthesis of feedback controllers that seek real- and reactive-power inverter setpoints corresponding to AC optimal power flow (OPF) solutions. The objective is to bridge the temporal gap between long-term system optimization and real-time inverter control, and enable seamless PV-owner participation without compromising system efficiency and stability. The design of the controllers is grounded on a dual $\epsilon$-subgradient method, and semidefinite programming relaxations are advocated to bypass the non-convexity of AC OPF formulations. Global convergence of inverter output powers is analytically established for diminishing stepsize rules for cases where: i) computational limits dictate  asynchronous updates of the controller signals, and ii) inverter reference inputs may be updated at a faster rate than the power-output settling time.
\end{abstract}
\begin{IEEEkeywords}
Distribution systems, photovoltaic inverter control, distributed optimization and control; optimal power flow. 
\end{IEEEkeywords}

\section{Introduction}
\label{sec:Introduction}

Present-generation residential photovoltaic (PV) inverters typically operate in a distributed and uncoordinated fashion, with the primary objective of maximizing the power extracted from PV arrays.  With the increased deployment of behind-the-meter PV systems,  an upgrade of medium- and low-voltage distribution-system operations and controls is required to address emerging efficiency, reliability, and power-quality concerns~\cite{Liu08,Woyte06}.  To this end, several architectural frameworks have been proposed for PV-dominant distribution systems to broaden the objectives of inverter real-time control, and enable inverters to partake in distribution-network optimization tasks~\cite{Aliprantis13, Tonkoski11, Farivar12,OID}. 

Past works have  addressed the design of \emph{distributed} real-time inverter-control strategies
to regulate the delivery of real and reactive power based on local measurements, so that terminal voltages are within
acceptable levels~\cite{Aliprantis13, Tonkoski11}. On a different time scale, centralized and distributed optimal power flow (OPF) algorithms have been proposed to compute optimal steady-state inverter setpoints, so that power losses and voltage deviations are minimized and economic benefits to end-users providing ancillary services are maximized~\cite{Khodr07, Paudyal11,Farivar12, OID_TEC,Dallanese-TSG13,Capitanescu14}. 

In an effort to bridge the temporal gap between real-time control and network-wide steady-state optimization, this paper addresses the synthesis of feedback controllers that seek optimal PV-inverter power setpoints corresponding to AC OPF solutions. 
The guiding motivation is to ensure that PV-system operation and control strategies are adaptable to changing ambient conditions and loads, and enable seamless end-user participation without compromising system efficiency. The proposed feedback controllers continuously pursue solutions of the formulated OPF problem by dynamically updating the setpoints based on current system outputs and problem parameters. This presents significant improvements over state-of-the-art distributed OPF approaches wherein reference signals for the inverters are updated only upon convergence of the distributed algorithm. In this setting, it is evident that if problem parameters or inputs change during the computation, broadcast, and implementation of the distributed solution of the OPF, the inverter would implement outdated setpoints.

Prior efforts in this direction include continuous-time feedback controllers that seek Karush-Kuhn-Tucker conditions for optimality  developed in~\cite{Jokic09}, and applied to solve an economic dispatch problem for bulk power systems in~\cite{Jokic_JEPES}.  Recently, modified automatic generation and frequency control methods that incorporate optimization objectives corresponding to DC OPF problems have been proposed for lossless bulk power systems in e.g.,~\cite{NaLi_ACC14,Chen_CDC14,Papachristodoulou14}. A heuristic based on saddle-point-flow methods is utilized in~\cite{Elia-Allerton13} to synthesize controllers seeking AC OPF solutions. Strategies that integrate economic optimization within droop control for islanded lossless microgrids are developed in~\cite{Dorfler14}. In a nutshell, these approaches are close in spirit to the seminal work~\cite{Arrowbook58}, where dynamical systems that serve as proxies for optimization variables and multipliers are synthesized to evolve in a continuous-time gradient-like fashion to the saddle points of the Lagrangian function associated with a convex optimization problem. For DC OPF, a heuristic comprising continuous-time dual ascent and discrete-time reference-signal updates is proposed in~\cite{Hirata15}; where, local stability of the resultant closed-loop system is also established.

Distinct from past efforts~\cite{Jokic_JEPES, Hirata15,NaLi_ACC14,Chen_CDC14, Papachristodoulou14, Dorfler14, Elia-Allerton13}, this work leverages dual $\epsilon$-subgradient methods~\cite{Larsson03,Kiwiel04}, to develop a feedback controller that steers the inverter output powers towards the solution of an AC OPF problem. A semidefinite programming (SDP) relaxation is advocated to bypass the non-convexity of the formulated AC OPF problem~\cite{Bai08,LavaeiLow, Dallanese-TSG13}. The proposed scheme involves the update of dual and primal variables in a discrete-time fashion, with the latter constituting the reference-input signals for the PV inverters. Convergence of PV-inverter-output powers to the solution of the formulated OPF problem is analytically established for settings where: \emph{i)} in an effort to bridge the time-scale separation between optimization and control, the reference inputs may be updated at a faster rate than the power-output settling time; and, \emph{ii)} due to inherent 
computational limits related to existing SDP solvers, the controller signals are updated asynchronously.  Although the present paper focuses on the case where an SDP relaxation is utilized to bypass the non-convexity of OPF tasks, the proposed synthesis procedure can be utilized to develop controllers that provably drive the inverter output to solutions of various convex relaxations~\cite{Farivar13} and linear approximations~\cite{Coffrin14,Bolognani15,DG-twrpd14} of the OPF problem. 

Overall, the proposed framework considerably broadens the approaches of~\cite{Jokic_JEPES, Hirata15,NaLi_ACC14,Chen_CDC14, Papachristodoulou14, Dorfler14, Elia-Allerton13} by: \emph{i)} considering AC OPF setups; \emph{ii)} incorporating PV-inverter operational constraints; \emph{iii)}~accounting for communication constraints which naturally lead to discrete-time controller updates; and, \emph{vi)} accounting for computational limits which involves an asynchronous update of the controller signals. It is also shown that the controller affords a distributed implementation, and requires limited message passing between the PV systems and the utility. 

The remainder of this paper is organized as follows. Section~\ref{sec:ProblemFormulation} outlines the problem formulation, while the PV controller is developed in Section~\ref{sec:Subgradient}. Section~\ref{sec:ApplicationExample} elaborates on the distributed implementation of the proposed control architecture. Numerical tests are reported in Section~\ref{sec:NumericalResults}, and conclusions are provided in Section~\ref{sec:Conclusions}.

\section{Problem Formulation}
\label{sec:ProblemFormulation}
Dynamical models and relevant  formulations for optimizing inverter setpoints are outlined for a general networked dynamical system in Section~\ref{sec:General}, and tailored to real-time PV-inverter control 
in Section~\ref{sec:Regulation}.\footnote{\emph{Notation}. Upper-case (lower-case) boldface
letters will be used for matrices (column vectors); $(\cdot)^\sfT$ for transposition; $(\cdot)^*$ complex-conjugate; and, $(\cdot)^\sfH$ complex-conjugate transposition;  $\Re\{\cdot\}$ and $\Im\{\cdot\}$ denote the real and imaginary parts of a complex number, respectively; $\mathrm{j} := \sqrt{-1}$. $\trace(\cdot)$ the matrix trace; $\mathrm{rank}(\cdot)$ the matrix rank; $|\cdot|$ denotes the magnitude of a number or the cardinality of a set; $\mathrm{vec}(\bX)$ returns a vector stacking the columns of matrix $\bX$, and $\mathrm{bdiag}(\{\bX_i\})$ forms a block-diagonal matrix. $\mathbb{R}^N$ and $\mathbb{C}^N$ denote the spaces of $N\times1$ real-valued and complex-valued vectors, respectively; $\mathbb{N}$ the set of natural numbers; and, $\mathbb{H}_+^{N \times N}$ denotes the space of $N \times N$ positive semidefinite Hermitian matrices. Given vector $\bx$ and square matrix $\bX$,  $\|\bx\|_2$ denotes the Euclidean norm of $\bx$, and $\|\bV\|_2$ the (induced) spectral norm of matrix $\bX$. $[\bx]_i$ ($[\bf(\bx)]_i$) points to the $i$-th element of a vector $\bx$ (vector-valued function $\bbf(\bx)$). $\dot{\bx}(t)$ is the time derivative of $\bx(t)$. Given a scalar function $f(\bx): \mathbb{R}^n \rightarrow \mathbb{R}$, $\nabla_\bx f(\bx)$ returns the gradient $[\frac{\partial f}{\partial x_1}, \ldots, \frac{\partial f}{\partial x_n}]^\sfT$. For a continuous function $f(t)$, $f[t_k]$ denotes its value sampled at $t_k$. Finally, $\bI_N$ denotes the $N \times N$ identity matrix; and, $\mathbf{0}_{M\times N}$, $\mathbf{1}_{M\times N}$ the $M \times N$ matrices with all zeroes and ones, respectively.}

\subsection{General problem setup}
\label{sec:General}

Consider $N_D$ dynamical systems described by
\begin{subequations}
\label{eq:sys1}
\begin{align}
\dot{\bx}_i(t) & = \bbf_i \Big( \bx_i(t), \bd_i(t), \bu_i(t) \Big) \label{eq:sys1-st} \\
\by_i(t) & =  \br_i \Big( \bx_i(t), \bd_i(t) \Big) , \quad i \in \mathcal{N}_D := \{1, \ldots, N_D \} \label{eq:sys1-obs} 
\end{align}
\end{subequations}
where: $\bx_i(t) \in \mathbb{R}^{n_{x,i}}$ is the state of the $i$-th dynamical system at time $t$; $\by_i(t) \in \cY_i \subset \mathbb{R}^{n_{y,i}}$ is the measurement of state $\bx_i(t)$ at time $t$; $\bu_i(t) \in \cY_i$ is the reference input; and $\bd_i(t) \in \cD_i  \subset \mathbb{R}^{n_{d,i}}$ is the exogenous input. Finally, $\bbf_i: \mathbb{R}^{n_{x,i}} \times \mathbb{R}^{n_{z,i}} \times \mathbb{R}^{n_{d,i}} \times \mathbb{R}^{n_{y,i}} \rightarrow \mathbb{R}^{n_{x,i}}$ and $\br_i:  \mathbb{R}^{n_{x,i}} \times \mathbb{R}^{n_{d,i}}  \rightarrow  \mathbb{R}^{n_{y,i}}$ are arbitrary (non)linear functions. The following system behavior for given finite exogenous inputs and reference signals is assumed.

\vspace{.1cm}

\begin{assumption}
\label{ass:DynSystems}
For given \emph{constant} exogenous inputs $\{\bd_i \in \cD_i\}_{i \in \cN_D}$ and reference signals $\{\bu_i \in \cY_i\}_{i \in \cN_D}$, there exist equilibrium points $\{\bx_i\}_{i \in \cN_D}$ for~\eqref{eq:sys1} that satisfy: 
\begin{subequations}
 \label{eq:sys1-equilibrium}
\begin{align}
\mathbf{0} & = \bbf_i \left(\bx_i, \bd_i, \bu_i \right) \label{eq:sys1-st-equilibrium} \\
\bu_i & =  \br_i \left( \bx_i, \bd_i \right), \quad i \in \cN_D  \, . \label{eq:sys1-obs-equilibrium} 
\end{align}
Notice that in~\eqref{eq:sys1-obs-equilibrium} the equilibrium output coincides with the commanded input $\bu_i$; that is, $\by_i = \bu_i$. 
These equilibrium points are locally asymptotically stable~\cite{KhalilBook}.  \hfill $\Box$ 
\end{subequations}
\end{assumption}

\vspace{.1cm}

For given exogenous inputs $\{\bd_i \in \cD_i\}_{i \in \cN_D}$, consider the following optimization problem: 
\begin{subequations}
\label{eq:p1}
\begin{align}
\mathrm{(P1)} \hspace{0cm}  \min_{\bV \in \cV, \{\bu_i \in \cY_i \}}& \,\, H(\bV) + \sum_{i \in \cN_D} G_i(\bu_i) \label{eq:p1-cost} \\
& \hspace{-1.8cm} \mathrm{subject~to}  \,\, \bh_i(\bV) - \bu_i + \bd_i  = \mathbf{0} , \, \forall \, i \in \cN_D \label{eq:p1-eq} 
\end{align}
\end{subequations}
where $\cV \subset \mathbb{H}_+^{n_V \times n_V}$  is a convex, closed, and bounded subset of the cone of positive semidefinite (Hermitian) matrices; function $H(\bV): \mathbb{H}_+^{n_V \times n_V} \rightarrow \mathbb{R}$ is known, convex and finite over $\cV$; $G_i(\bu_i): \mathbb{R}^{n_{y,i}} \rightarrow \mathbb{R}$ is strongly convex and finite over $\cY_i$; and, the vector-valued function $\bh_i(\bV): \mathbb{H}_+^{n_V \times n_V} \rightarrow \mathbb{R}^{n_{y,i}}$ is affine. Finally, sets $\{\cY_i\}_{i \in \cN_D}$, which define the space of possible reference inputs for the dynamical systems, are assumed to comply to the following requirement.

 \vspace{.1cm}
\noindent   
\begin{assumption}
\label{ass:compactSetY}
Sets $\{\cY_i\}_{i \in \cN_D}$ are convex and compact. Further,  $\mathrm{(P1)}$ has a non-empty feasible set and a finite optimal cost. \hfill $\Box$ 
\end{assumption}

\vspace{.1cm}

\noindent With these assumptions, problem $\mathrm{(P1)}$ is a \emph{convex} program; moreover, it can be reformulated into a standard SDP form by resorting to the epigraph form of the cost function. 

It is evident from~\eqref{eq:sys1-obs-equilibrium} that $\mathrm{(P1)}$ defines the optimal operating setpoints of the dynamical systems~\eqref{eq:sys1} in terms of steady-state outputs~\cite{Jokic09,Hirata15}. In fact, by utilizing the optimal solution $\{\bu_i^{\textrm{opt}}\}_{i \in \cN_D}$ of $\mathrm{(P1)}$ as reference inputs, it follows from~\eqref{eq:sys1-obs-equilibrium} that each system output will eventually be driven to the point $\by_i = \bu_i^{\textrm{opt}}$. Function~\eqref{eq:p1-cost} captures costs incurred by the steady-state outputs, as well as costs associated with matrix variable $\bV$, which couples the steady-state system outputs $\{\by_i = \bu_i\}_{i \in \cN_D}$ through the linear equality constraints~\eqref{eq:p1-eq}.

In principle, $\mathrm{(P1)}$ could be solved centrally by a system-level control unit, which subsequently dispatches the reference signals $\{\bu_i^{\textrm{opt}}\}_{i \in \cN_D}$ for the dynamical systems. In lieu of a centralized solution of $\mathrm{(P1)}$, the \emph{objective} here is to design a \emph{distributed} feedback controller for the dynamical systems~\eqref{eq:sys1}, so that the resultant closed-loop system is globally convergent to an equilibrium point $\{\bx_i\}_{i \in \cN_D}$, $\{\by_i =  \br_i ( \bx_i, \bd_i )\}_{i \in \cN_D}$, where the values $\{\by_i\}_{i \in \cN_D}$ of the steady-state outputs coincide with the optimal solution $\{\bu_i^{\textrm{opt}}\}_{i \in \cN_D}$ of $\mathrm{(P1)}$.

\subsection{PV-inverter output regulation to OPF solutions}
\label{sec:Regulation}

The task of regulating the power output of PV inverters is outlined in this section, and cast within the framework of~\eqref{eq:sys1}-\eqref{eq:p1}. In this regard,~\eqref{eq:sys1}-\eqref{eq:sys1-equilibrium} will model the inverter dynamics~\cite[Ch.~8]{Iravanibook10},~\cite{Irminger12}; while  OPF  will be formulated in the form~\eqref{eq:p1} by leveraging  SDP relaxation techniques~\cite{Bai08,LavaeiLow}.

\emph{Network.}  Consider a distribution system comprising $N+1$ nodes collected in the
set $\cN$, and lines represented by the set of undirected edges $\cE := \{(m,n): m, n \in \cN\}$. The set $\cN := \{0, 1, \ldots, N\}$ is partitioned as $\cN = \{0\} \cup \cN_D \cup \cN_O$, where: node $0$ denotes the secondary of the step-down transformer; inverter-interfaced PV systems are located at nodes $\cN_D = \{1, \ldots, N_D\}$ [cf.~\eqref{eq:sys1}]; and, $\cN_O := \{N_D+1, \ldots, N\}$ collects nodes with no power generation. For simplicity of exposition, the framework is outlined for a balanced system; however, the proposed framework can be extended to unbalanced multi-phase systems as explained in Appendix~\ref{sec:three-phase}. 

\vspace{.2cm}

\emph{Dynamics of PV inverters.} Equation~\eqref{eq:sys1-st} is utilized to model the dynamics of PV inverters, regulating real- and reactive output powers to prescribed setpoints. For example, relevant dynamical models for inverters operating in a grid-connected mode are discussed in e.g.,~\cite[Ch.~8]{Iravanibook10} and~\cite{Irminger12}. These models can be conveniently cast within~\eqref{eq:sys1}-\eqref{eq:sys1-equilibrium} as shown next. 

\noindent  $\bullet$ Let $p_i(t) := E_i(t) \cos(\omega t + \phi_i(t)) i_i(t)$ and $q_i(t) := E_i(t) \cos(\omega t + \phi_i(t) - \pi/2) i_i(t)$ denote the instantaneous output real and reactive powers of inverter $i \in \cN_D$, respectively, where $\omega$ is the grid frequency, $v_i(t) := E_i(t) \cos(\omega t + \phi_i(t))$ the voltage waveform, and $i_i(t)$ is the  current injected. Further, let $P_i(t)$ and $Q_i(t)$ denote averages of the instantaneous output real and reactive powers over an AC cycle; that is,  
\begin{align}
P_i(t) := \frac{\omega}{2 \pi} \int_{t- \frac{\omega}{2 \pi}}^t \hspace{-.3cm} p_i(\tau) \mathrm{d} \tau , \,\,\,\,\,
Q_i(t)  := \frac{\omega}{2 \pi} \int_{t- \frac{\omega}{2 \pi}}^t \hspace{-.3cm}  q_i(\tau) \mathrm{d}  \tau . \label{eq:PQcycle} 
\end{align} 
Then, the \emph{state} of system~\eqref{eq:sys1} is $\bx_i(t)  := [P_{i}(t),Q_{i}(t)]^\sfT$.  

\noindent  $\bullet$  Vector ${\bf u}_i (t) = {\bf u}_i$ collects the constant \emph{commanded} real and reactive powers for inverter $i$. By~\eqref{eq:sys1-equilibrium}, inverters regulate the output powers to the commanded setpoints $\bu_i$; see e.g.,~\cite[Ch.~8]{Iravanibook10},~\cite{Irminger12}.

\noindent  $\bullet$ Let ${P}_{\ell,i}(t)$ and ${Q}_{\ell,i}(t)$ denote the demanded real and reactive
 loads at node $i \in \cN$. Then, vector $\bd_i(t)$ is set to be $\bd_i(t) :=  [{P}_{\ell, i}(t),{Q}_{\ell, i}(t)]^\sfT$ for all $i \in \cN \backslash \{0\}$.   

\noindent  $\bullet$ By setting 
\vspace{-.3cm}
\begin{align}
\br_i ( \bx_i(t), \bd_i(t) ) = \bx_i(t)
\end{align}
~\eqref{eq:sys1-obs} equates the state with the measurement of the inverter output powers.

\vspace{.2cm}

\emph{Steady-state OPF problem}. Let $V_i := (E_i/\sqrt{2}) e^{\imj \phi_{i}} \in \mathbb{C}$ be the phasor representation of 
the steady-state voltage at node $i \in \cN$. Similarly, let $I_i \in \mathbb{C}$ denote the phasor for the current injected at node $i \in \cN$, and define $\bi := [I_0, \ldots, I_N]^\sfT \in \mathbb{C}^{N+1}$ and $\bv := [V_0, \ldots, V_N]^\sfT \in
\mathbb{C}^{N+1}$. Then, using Ohm's and Kirchhoff's circuit laws, the linear
relationship $\bi = \bY \bv$ can be established, where $\bY \in \mathbb{C}^{N+1 \times N+1}$ is the admittance matrix formed
based on the distribution-network topology and the $\pi$-equivalent circuits 
of lines $\cE$~\cite[Ch. 6]{Kerstingbook}; see e.g.,~\cite{Bai08,LavaeiLow,Dallanese-TSG13,Robbins15} for details on the construction of matrix $\bY$.

For prevailing ambient conditions, let $P_i^{\textrm{av}} \geq 0$ denote the \emph{available} real power for 
the inverter $i \in \cN_D$. The available power is a function of the incident irradiance, and corresponds to the maximum power point of the PV array. When PV-systems operate at unity power factor~\cite{OID},  a set of challenges related to power quality and reliability in distribution systems may emerge for sufficiently high levels if deployed PV capacity~\cite{Liu08}. For instance, overvoltages may be experienced during periods when PV generation exceeds the household demand~\cite{Liu08}, while  fast-variations in the PV-output tend to cause transients that lead to wear-out of legacy switchgear~\cite{Woyte06}.  Efforts to ensure reliable operation of existing distribution systems with increased behind-the-meter PV generation are  focused on the possibility of inverters providing reactive power compensation and/or curtailing real power~\cite{Aliprantis13,Tonkoski11,Farivar12,OID}. In the most general setting, the set of operating points for PV inverters providing ancillary services can be specified as:
\begin{align}
& \cY_i =  \left\{({P}_{i}, {Q}_{i} ) \hspace{-.1cm} :  {P}_{i}^{\textrm{min}}  \leq {P}_{i}  \leq  P_{i}^{\textrm{av}}, {Q}_{i}^2  \leq  S_{i}^2 - {P}_{i}^2, \right. \nonumber \\
& \hspace{4cm} \left. \mathrm{~and~} | {Q}_{i} | \leq (\tan \theta) {P}_{i} \label{mg-PV} 
 \right\}
\end{align}
where $S_i$  is the rated apparent power, and the last inequality is utilized to enforce a minimum power factor of $\cos \theta$. Parameters $\theta$ and ${P}_{i}^{\textrm{min}}$ can be conveniently tuned to account for the following strategies: 

\noindent $\mathrm{(c1)}$ Reactive power compensation: ${P}_{i}^{\textrm{min}} = P_{i}^{\textrm{av}}$, $\theta \in (0, \pi/2]$;

\noindent $\mathrm{(c2)}$ Active power curtailment:  ${P}_{i}^{\textrm{min}} \in [0, P_{i}^{\textrm{av}})$, $\theta = 0$; and,

\noindent $\mathrm{(c3)}$ Active and reactive control: ${P}_{i}^{\textrm{min}} \in [0, P_{i}^{\textrm{av}}]$, $\theta \in (0, \pi/2]$.

\begin{figure}[t]
\begin{center}
\subfigure[]{\includegraphics[width=2.1cm]{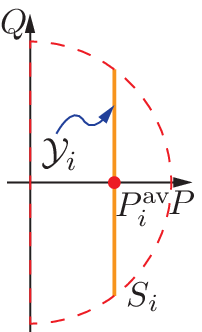}} \hspace{.5cm}
\subfigure[]{\includegraphics[width=2.1cm]{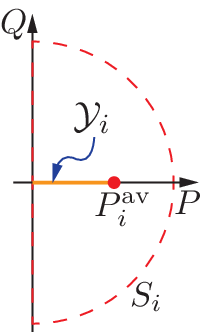}} \hspace{.5cm}
\subfigure[]{\includegraphics[width=2.1cm]{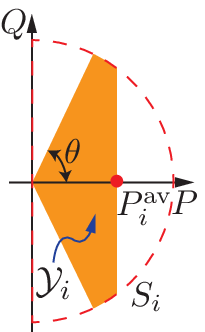}}
\end{center}
\vspace{-.4cm}
\caption{Operating regions $\cY_i$ for PV inverter under: (a) reactive power compensation~\cite{Aliprantis13}; (b) real power curtailment~\cite{Tonkoski11}; and, (c) combined real- and reactive-power control~\cite{OID}.}
\label{Fig:OIDregions}
\vspace{-.4cm}
\end{figure}

\noindent The PV-inverter operating regions involved by strategies $\mathrm{(c1)}$--$\mathrm{(c3)}$ are illustrated in Figure~\ref{Fig:OIDregions}. It is evident that sets $\{\cY_i\}$ adhere to \emph{Assumption~\ref{ass:compactSetY}}.

For given load and ambient conditions, a prototypical OPF formulation for optimizing the
steady-state operation of a distribution system  is given as follows:
\begin{subequations} 
\label{Pmg}
\begin{align} 
 \mathrm{(OPF)}  \hspace{1.8cm} & \hspace{-1.7cm} \min_{\bv, \bi,  \{P_i, Q_i \} } \,\, H(\bv) +\sum_{i \in \cN_D} G_i(P_i, Q_i) \label{mg-cost} \\ 
\mathrm{subject\,to} \,\, & \bi = \bY \bv,  \mathrm{and}  \nonumber  \\ 
& \hspace{-1.3cm} V_i I_i^* = P_i - P_{\ell,i} + \mathrm{j} (Q_i - Q_{\ell,i}), \hspace{.1cm}  \forall \, i \in \cN_D \label{mg-balance-I} \\
& \hspace{-1.35cm} V_n I_n^* = - P_{\ell,n} - \mathrm{j} Q_{\ell,n}, \hspace{1.55cm}  \forall \, n \in \cN_O \label{mg-balance-L} \\
& \hspace{-1.35cm} V^{\mathrm{min}} \leq |V_i| \leq V^{\mathrm{max}}  \hspace{1.95cm}  \forall \, i \in \cN  \hspace{-.2cm} \label{mg-Vlimits} \\
& \hspace{-.9cm} \bu_i \in \cY_i \hspace{3.45cm} \forall \, i \in \cN_D \label{mg-PVp} 
\end{align}
\end{subequations}
where $V^{\mathrm{min}}$ and $V^{\mathrm{max}}$ are prescribed voltage limits (e.g., ANSI C84.1 limits); the constraint on $|V_0|$ is left implicit;~\eqref{mg-PVp} specifies the feasible inverter operating region [cf. Figure~\ref{Fig:OIDregions}];  and, equalities~\eqref{mg-balance-I}--\eqref{mg-balance-L} capture the power-balance equations for nodes with and without inverters, respectively. For nodes without loads (e.g., utility poles), one clearly has that ${P}_{\ell,i} = {Q}_{\ell,i} = 0$. 

Function $H(\bv)$ can capture various network-oriented performance objectives that the distribution system operator (DSO) may pursue. For example, the DSO may aim to minimize the power losses on the distribution lines, voltage magnitude deviations from nominal, and/or the power drawn from the substation~\cite{Khodr07, Paudyal11}. On the other hand, function $G_i(P_i, Q_i)$ models PV-inverter costs/rewards for ancillary service provisioning such as real power curtailment and/or reactive power compensation~\cite{Farivar12,OID,Dorfler14}; for example, this function can be set to $G_i(P_i, Q_i) = a_i (P_{i}^{\textrm{av}} - P_i)^2 + b_i (P_{i}^{\textrm{av}} - P_i) + c_i Q_i^2 + d_i |Q_i|$, with $a_i, b_i, c_i, d_i$ denoting market-oriented coefficients, to maximize the amount of power provided by PV systems.  Finally, notice that additional constraints such as thermal limits may be naturally accommodated in~\eqref{Pmg}.

It is well-known that the OPF problem~\eqref{Pmg} is \emph{nonconvex}, and thus hard to solve to global optimality in both centralized and distributed setups. Further, given that the problem is nonconvex, convergence of distributed algorithms (derived, e.g., via Lagrangian decomposition techniques) is not always guaranteed and needs to be  established.  Since the objective  of this work is to develop \emph{distributed} controllers so that inverter output powers are \emph{provably} convergent to OPF solutions, a convex reformulation of the OPF task is considered next.

\vspace{.2cm}

\emph{SDP relaxation of the OPF problem}. To formulate an SDP relaxation of the prototypical steady-state OPF problem~\eqref{Pmg}, consider expressing steady-state powers and voltage magnitudes as linear functions of the
outer-product matrix $\bV := \bv \bv^\sfH$~\cite{Bai08,LavaeiLow}, and define matrix 
$\mathbf{Y}_i := \be_i \be_i^\sfT \mathbf{Y}$ per node $i$, where $\{\mathbf{e}_{i}\}_{i \in \cN}$
denotes the canonical basis of $\mathbb{R}^{|\cN|}$. Using
$\mathbf{Y}_i$, form the Hermitian matrices $\bPhi_{i} :=
\frac{1}{2} (\bY_i + \bY_i^\sfH)$, $\bPsi_{i} := \frac{\textrm{j}}{2} (\bY_i -
\bY_i^\sfH) $, and $\bUpsilon_{i} := \be_i \be_i^\sfT$. Then, the balance equations for real and reactive powers at node 
$i \in \cN_D$ can be expressed as $\trace(\bPhi_i \bV)   =  {P}_{i} - {P}_{\ell,i}$ and 
$\trace(\bPsi_i \bV) = {Q}_{i} - {Q}_{\ell,i}$, respectively. To reformulate the OPF in the form~\eqref{eq:p1}, consider setting $\bu_i  = [{P}_{i},{Q}_{i}]^\sfT $, $\bd_i  = [{P}_{\ell, i},{Q}_{i}]^\sfT$, and
\begin{align}
\label{eq:correspondence_h}
\bh_i(\bV)  = [\trace(\bPhi_i \bV), \trace(\bPsi_i \bV)]^\sfT . \hspace{-.1cm}
\end{align} 
Additionally, define the following convex set: 
\begin{align}
& \cV :=  \{\bV: \bV \succeq \mathbf{0}, V_{\mathrm{min}}^2 \leq \trace(\bUpsilon_i \bV)  \leq V_{\mathrm{max}}^2 \,\forall \, i \in \cN  \, \nonumber  \\
& \mathrm{~and~} \trace(\bPhi_i \bV)   =  - {P}_{\ell,i}, \trace(\bPsi_i \bV)   =  - {Q}_{\ell,i} \, \forall i \in \cN_O \} \, .  \label{eq:V_region}  
\end{align} 
With these definitions, problem~\eqref{Pmg} can be \emph{equivalently} expressed as follows
\begin{subequations}
\label{eq:opf2}
\begin{align}
 \hspace{0cm}  \min_{\bV \in \cV, \{\bu_i \in \cY_i \}}& \,\, H(\bV) + \sum_{i \in \cN_D} G_i(\bu_i) \label{eq:opf2-cost} \\
& \hspace{-1.5cm} \mathrm{subject~to}  \,\, \bh_i(\bV) - \bu_i + \bd_i  = \mathbf{0} , \, \forall \, i \in \cN_D \label{eq:opf2-eq} \\
& \hspace{.2cm}  \mathrm{rank}(\bV) = 1. \label{eq:opf2-rank}
\end{align}
\end{subequations}
On par with~\eqref{Pmg}, problem~\eqref{eq:opf2} is nonconvex because of the rank constraint; however in the spirit of the SDP relaxation, the constraint~\eqref{eq:opf2-rank} can be dropped. Notice that, once the constraint~\eqref{eq:opf2-rank} is dropped, the resultant SDP relaxation of the OPF problem is in the form~\eqref{eq:p1}. If the optimal matrix $\bV^{\mathrm{opt}}$ of the relaxed problem has $\mathrm{rank}(\bV^{\mathrm{opt}}) = 1$, then the resultant power flows are globally optimal~\cite{Bai08, LavaeiLow}. Sufficient conditions for this relaxation to be exact for radial and balanced systems are provided in~\cite{Zhang14}, while its applicability to unbalanced multiphase systems is investigated in~\cite{Dallanese-TSG13}. 

In this setup, the \emph{objective} of the feedback controller that will be designed in the following section, is to drive the inverter outputs  $\{\by_i(t) = [P_{i}(t),Q_{i}(t)]^\sfT\}_{i \in \cN_D}$ to the optimal solution $\{\bu_i^{\textrm{opt}}\}_{i \in \cN_D}$ of the OPF problem.

\section{Feedback Controller}
\label{sec:Subgradient}

Dual $\epsilon$-subgradient methods are leveraged in Section~\ref{sec:controllerSynthesis} to 
synthesize controllers for systems~\eqref{eq:sys1} whose outputs track recursive solvers of $\mathrm{(P1)}$. Applications to the real-time PV-inverter control problem are discussed in Section~\ref{sec:ApplicationExample}. 

To streamline proofs of relevant analytical results, it will be convenient to express the linear equality constraints~\eqref{eq:p1-eq} in a compact form. To this end, define $\bu := [\bu_1^\sfT, \ldots, \bu_{N_D}^\sfT]^\sfT$,  $\bd := [\bd_1^\sfT, \ldots, \bd_{N_D}^\sfT]^\sfT$ and $\bh(\bV) := [\bh_1^\sfT(\bV), \ldots, \bh_{N_D}^\sfT(\bV)]^\sfT$. Then, constraints~\eqref{eq:p1-eq} can be compactly expressed as $\bh(\bV) = \bu - \bd$. 

\subsection{Primer on dual gradient method}
\label{sec:dualgradient}

Consider the Lagrangian corresponding to~\eqref{eq:p1}, namely:  
\begin{align}
L\left(\bV, \{\bu_i\}, \{\blambda_i\}\right) &:=  \,H(\bV) + \sum_{i \in \cN_D} G_i(\bu_i) \nonumber \\
& + \sum_{i \in \cN_D} \blambda_i^\sfT \left( \bh_i(\bV) - \bu_i + \bd_i \right) \label{eq:lagrangian}
\end{align}
where $\blambda_i \in \mathbb{R}^{n_{y,i}}$ denotes the Lagrange multiplier associated with~\eqref{eq:p1-eq}. Based on~\eqref{eq:lagrangian}, the dual function and the dual problem are defined as follows (see, e.g.,~\cite{Bertsekas_ConvexAnalysis}): 
\begin{align}
q(\{\blambda_i\}) := \min_{\bV \in \cV, \{\bu_i \in \cY_i\}_{i \in \cN_D}} L(\bV, \{\bu_i\}, \{\blambda_i\})  \label{eq:dualFunction}
\end{align}
\begin{align}
q^{\mathrm{opt}} := \max_{\{\blambda_i\}_{i \in \cN_D}} \,\, q(\{\blambda_i\})   \, .\label{eq:optimalDualFunction}
\end{align}

Regarding the optimal Lagrange multipliers, the following technical requirement is presumed in order to guarantee their existence and uniqueness; see e.g.,~\cite{Wachsmuth13}. 

\vspace{.1cm}

\begin{assumption}
\label{ass:constraintQualification}
Vectors
\begin{align}
\nabla_{[\textrm{vec}^\sfT (\bV), \bu^\sfT ]^\sfT} [\bh(\bV) + \bg(\bu,\bd)]_i , \,\, i = 1, \ldots, \sum_{i \in \cN_D} n_{y,i} \label{eq:cq}
\end{align}
are linearly independent. \hfill $\Box$ 
\end{assumption}

\vspace{.1cm}

\noindent Section~\ref{sec:ApplicationExample} will elaborate on how condition~\eqref{eq:cq} can be checked in the OPF context.  Under current modeling assumptions, it follows that the duality gap is zero, and the dual function $q(\{\blambda_i\})$ is concave, differentiable, and it has a continuous first derivative~\cite{Cheng87}. Consider then  utilizing a gradient method to solve the dual problem, which amounts to iteratively performing~\cite{Cheng87}:  
\begin{subequations} 
\label{eq:dualsubgradient_basic}
\begin{align}
\hspace{-.2cm} \{\bV[k], \{\bu_i[k]\}_{i \in \cN_D}\} & \nonumber \\
& \hspace{-1.55cm} = \arg \min_{\bV \in \cV, \{\bu_i \in \cY_i\}} \, L(\bV, \{\bu_i\}, \{\blambda_i[k]\}) \label{eq:primal_basic} \\
&\hspace{-2.9cm}  \blambda_i[k+1]  = \blambda_i[k] + \alpha_{k+1} \nabla_{\blambda_i}  L(\bV[k], \{\bu_i[k]\}, \{\blambda_i\})\,  \label{eq:dual_ascent_basic}
\end{align}
\end{subequations}
where $k \in \mathbb{N}$ denotes the iteration index,  $\alpha_{k+1} \geq 0$ is the stepsize, and~\eqref{eq:dual_ascent_basic} is repeated for all $i \in \cN$. In particular, a non-summable but square-summable stepsize sequence is adopted in this paper~\cite{Kiwiel04}; that is, there exist sequences $\{\gamma_k\}_{k \geq 0}$ and $\{\eta_k\}_{k \geq 0}$ such that: \\
\noindent $\mathrm{(s1)}$ $\gamma_k \rightarrow 0$ as $k \rightarrow + \infty$, and $\sum_{k = 0}^{+ \infty} \gamma_k = + \infty$;  \\
\noindent $\mathrm{(s2)}$ $\gamma_k \leq \alpha_k \leq \eta_k$ for all $k \geq 0$; and, \\
\noindent $\mathrm{(s3)}$ $\eta_k \downarrow 0$ as $k \rightarrow + \infty$, and $\sum_{k = 0}^{+ \infty} \eta_k^2 < + \infty$.  \\
\noindent At iteration $k$, the same step-size $\alpha_k$ is utilized for all $i \in \cN$.  Exploiting the decomposablility of the Lagrangian, steps~\eqref{eq:dualsubgradient_basic}  can be equivalently expressed as:  
\begin{subequations}
\label{eq:dualsubgradient}
\begin{align}
& \hspace{-.2cm} \blambda_i[k+1] = \blambda_i[k] + \alpha_{k+1} \left( \bh_i(\bV[k]) - \bu_i[k] + \bd_i \right)  \label{eq:dual_ascent} \\
& \hspace{-.2cm} \bu_i [k+1] = \arg \min_{\bu_i \in \cY_i}  G_i(\bu_i) - \blambda_i^\sfT[k+1] \bu_i \label{eq:primal_y_2}  \hspace{-.2cm}  \\
& \hspace{-.2cm} \bV[k+1] = \arg \min_{\bV \in \cV} \, H(\bV) +  \sum_{i \in \cN_D} \blambda_i^\sfT[k+1] \, \bh_i(\bV) \label{eq:primal_V} 
\end{align}
\end{subequations} 
with~\eqref{eq:primal_y_2}--\eqref{eq:dual_ascent} performed for all $i \in \cN_D$. Finally, notice that from the compactness of sets $\cV$ and $\{\cY_i\}_{i = 1}^N$, it follows that there exists a scalar $G$, $0 \leq G < + \infty$, such that
\begin{align}
\| \bh(\bV[k]) - \bu[k] + \bd) \|_2 \leq G \, , \quad \forall \,\,\, k \in \mathbb{N} \, . \label{eq:bounded_gradient}
\end{align}

Using~\eqref{eq:bounded_gradient}, and a stepsize sequence $\{\alpha_k\}_{k\geq 0}$ satisfying $\mathrm{(s1)}$--$\mathrm{(s3)}$, it turns out that the dual iterates $\blambda_i[k]$ converge to the optimal solution $\blambda_i^{\mathrm{opt}}$ of the dual problem~\eqref{eq:optimalDualFunction}; that is, $\|\blambda^{\mathrm{opt}} - \blambda[k]\|_2 \rightarrow 0$ as $k \rightarrow \infty$~\cite[Prop.~8.2.6]{Bertsekas_ConvexAnalysis},~\cite{Cheng87,Kiwiel04}. Iterates $\bV[k]$ and $\{\bu_i[k]\}_{i \in \cN_D}$ become asymptotically feasible and their optimal values, $\bV^{\mathrm{opt}}$ and $\{\bu_i^{\mathrm{opt}}\}_{i \in \cN_D}$, can be recovered from~\eqref{eq:primal_V} and~\eqref{eq:primal_y_2}, respectively, once $\{\blambda_i^{\mathrm{opt}}\}_{i \in \cN_D}$ becomes available.  

Steps similar to~\eqref{eq:dualsubgradient} are typically adopted to enable a distributed solution of the OPF~\cite{Tse12,OID_TEC, Robbins15,Baldick99,Hug09,Erseghe14} as well as other resource allocation tasks such as the economic dispatch problem and residential load control~\cite{GatsisTSG12}. As illustrated in Figure~\ref{Fig:diagram}(a) and explained in detail in Section~\ref{sec:ApplicationExample}, updates~\eqref{eq:dual_ascent}-\eqref{eq:primal_y_2} are implemented \emph{at each individual PV system} while~\eqref{eq:primal_V} is performed \emph{at the DSO}. However, in conventional approaches, the optimal reference signals $\{\bu_i^{\textrm{opt}}\}_{i \in \cN_D}$ are implemented at the PV-inverters \emph{only} when the distributed algorithm converges to the optimal solution. It is evident that under this operating paradigm the  optimization and local control tasks operate \emph{at two different time scales}, with reference signals updated every time that the OPF problem is solved and implemented only when the inverter dynamics are in steady state. This motivates the development of control schemes that continuously pursue solutions of the OPF problem by dynamically updating the setpoints, based on current system outputs and problem parameters. This is described next.

\subsection{Controller synthesis}
\label{sec:controllerSynthesis}

Consider updates  performed at discrete time instants $t \in \{t_k, k \in \mathbb{N}\}$, with $\bV[t_k]$, $\{\bu_i[t_k]\}_{i \in \cN_D}$, and let $\{ \blambda_i[t_{k}]\}_{i \in \cN_D}$ denote the values of primal and dual variables, respectively, at time $t_k$. The following method accounts for the system dynamics in~\eqref{eq:sys1} while solving $(\mathrm{P1})$ with dual-gradient-based approaches. 

At time $t_{k}$, the system outputs are sampled as:
\begin{subequations}
\label{eq:dualsubgradient_system}
\begin{align}
\by_i[t_{k}] & =  \br_i \Big( \bx_i(t_{k}), \bd_i \Big) \, , \forall \, i \in \cN_D \label{eq:sys1-obs_sys}  
 \end{align}
and, they are utilized to update the dual variables as follows:    
\begin{align}
 \blambda_i[t_{k+1}]  & = \blambda_i[t_{k}] \nonumber \\
 & \hspace{-1.1cm} + \alpha_{k+1} \Big( \bh_i(\bV[t_{k}]) -  \by_i[t_{k}] + \bd_i  \Big)  \,   , \, \forall i  \in \cN_D. \label{eq:dual_ascent_sys}
\end{align}
Given $\blambda_i[t_{k+1}]$, the primal variables $\bV[t_{k+1}]$ and $\{ \bu_i[t_{k+1}]\}_{i \in \cN_D}$ are then updated as: 
\begin{align}
& \hspace{-.1cm} \bu_i [t_{k+1}] = \arg \min_{\bu_i \in \cY_i} G_i(\bu_i) - \blambda_i^\sfT[t_{k+1}] \bu_i \, .\label{eq:primal_y_sys}  \\
& \hspace{-.1cm} \bV[t_{k+1}]  = \arg \min_{\bV \in \cV} \, H(\bV) +  \sum_{i \in \cN_D}  \blambda_i^\sfT[t_{k+1}] \, \bh_i(\bV) \label{eq:primal_V_sys} 
\end{align}
\end{subequations}
Once~\eqref{eq:primal_y_sys} is solved, the vector-valued reference signal $\bu_i[t_{k+1}]$ is applied to the dynamical system~\eqref{eq:sys1-st} over the interval $(t_{k}, t_{k+1}]$; that is, $\bu_i(t) = \bu_i[t_{k+1}], t \in (t_{k}, t_{k+1}]$. At time $t_{k+1}$ the system outputs $\{\by_i[t_{k+1}] \}_{i \in \cN_D}$ are sampled again, and~\eqref{eq:dual_ascent_sys}--\eqref{eq:primal_y_sys} are repeated. Notice that, differently from standard dual gradient methods, variable $\bu_i[t_{k}]$ is replaced by the sampled system output $\by_i[t_{k}]$ in the ascent step~\eqref{eq:dual_ascent_sys}. 

\begin{figure*}[t] 
\begin{center}
\subfigure[]{\includegraphics[width=12.7cm ]{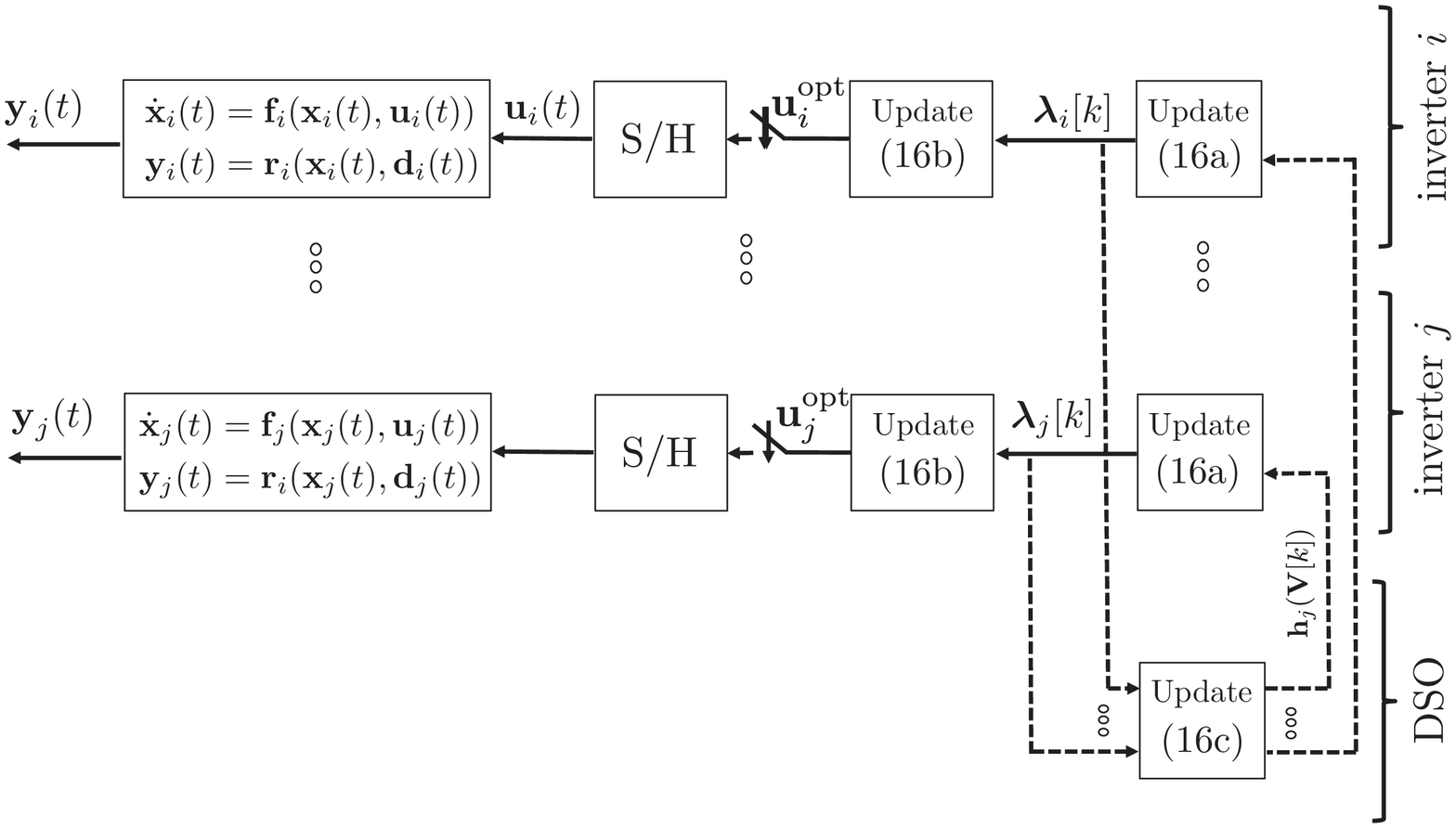}}
\subfigure[]{\hspace{-.45cm} \includegraphics[width=12.7cm ]{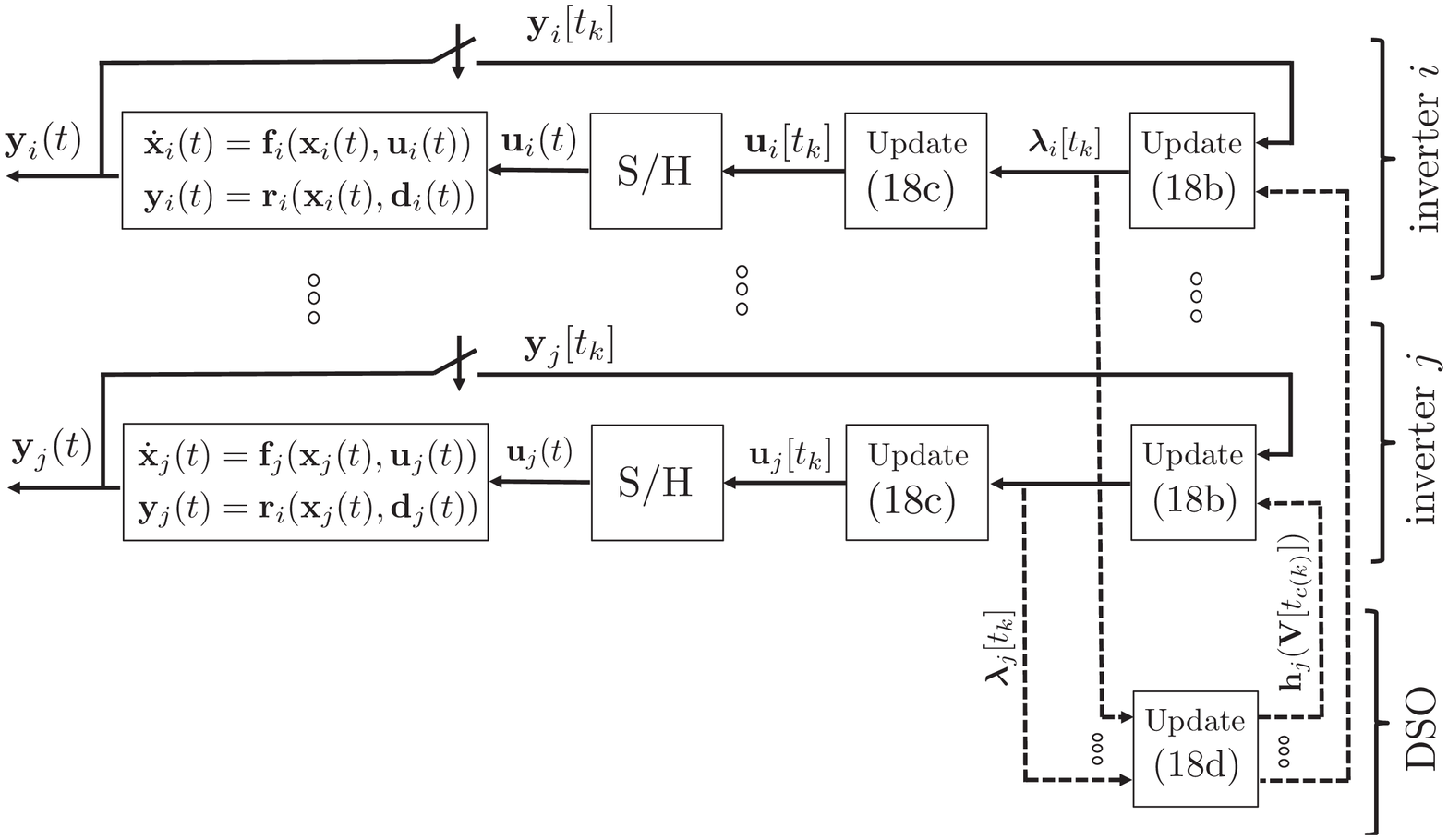}}
\end{center}
\vspace{-.3cm}
\caption{(a)~\emph{Conventional distributed optimization setup}: Problem~\eqref{eq:p1} is solved in a distributed fashion by using steps~\eqref{eq:dualsubgradient}; once the problem is solved (i.e., iterates in~\eqref{eq:dualsubgradient} have converged to the optimal primal and dual values), the optimal reference signals $\{\bu_i^{\textrm{opt}}\}_{i \in \cN_D}$ are dispatched to the PV-inverters. (b)~\emph{Proposed optimization-centric control architecture}: The discrete-time control signal $\bu_i[t_k]$ generated by the dynamic controller $i \in \cN_D$ is dynamically applied as an input to the inverter by utilizing a sample-and-hold (S/H) unit; the instantaneous inverter output is sampled and utilized for updating the control signals. The same architecture is utilized in the asynchronous case~\eqref{eq:dualsubgradient_system_as}, upon substituting~\eqref{eq:primal_V_sys} with~\eqref{eq:primal_V_sys_as}. As claimed in Theorem~1, the inverter outputs $\{\by_i(t)\}_{i \in \cN_D}$ converge to the solution of the OPF problem. Details on the distributed implementation are provided in Section~\ref{sec:ApplicationExample}.}
\label{Fig:diagram}
\vspace{-.3cm}
\end{figure*}

Steps~\eqref{eq:dual_ascent_sys}--\eqref{eq:primal_y_sys} in effect constitute the controller for the dynamical systems~\eqref{eq:sys1}. Specifically, the (continuous-time) reference signals $\{\bu_i(t)\}_{i \in \cN_D}$ produced by the controller have step changes at instants $\{t_k, k \in \mathbb{N}\}$, are left-continuous functions, and take the constant values $\{\bu_i[t_{k+1}]\}_{i \in \cN_D}$ over the time interval $(t_k, t_{k+1}]$. It is evident that if $\bu_i[t_k]$ converges to $\bu_i^{\textrm{opt}}$ as $k \rightarrow \infty$ (and thus $\bu_i(t) \rightarrow \bu_i^{\textrm{opt}}$ as $t \rightarrow \infty$), then $\by_i(t) \rightarrow \bu_i^{\textrm{opt}}$ as $t \rightarrow \infty$ by virtue of~\eqref{eq:sys1-equilibrium}. 

Suppose for now that the interval $(t_{k-1}, t_{k}]$ is large enough to allow the outputs $\{\by_i(t)\}_{i \in \cN_D}$ to converge to the commanded input $\{\bu_i[t_k]\}_{i \in \cN_D}$ [cf.~\eqref{eq:sys1-equilibrium}]. Under this ideal setup with a pronounced and tangible \emph{time-scale separation} between controller and system dynamics, one has that $\lim_{t \rightarrow t_{k}^-} \|\by_i(t) - \bu_i[t_k]\| = 0$, for all $k$ [cf.~\eqref{eq:sys1-equilibrium}], and  step~\eqref{eq:dual_ascent_sys} is replaced by $\blambda_i[t_{k+1}] = \blambda_i[t_{k}] + \alpha_{k+1} ( \bh_i(\bV[t_{k}]) - \bu_i[t_{k}] + \bd_i)$. Thus,~\eqref{eq:dualsubgradient_system} coincides with standard dual gradient method in~\eqref{eq:dualsubgradient}, and the convergence results in~\cite[Prop.~8.2.6]{Bertsekas_ConvexAnalysis},~\cite{Cheng87} carry over to this ideal setup. In this work, convergence of the system outputs $\{\by_i(t)\}_{i \in \cN_D}$ to the solution of $\mathrm{(P1)}$ is assessed in the more general case where update of reference signals may be performed \emph{faster than the systems' settling times} and \emph{asynchronously}, in order to achieve the following operational goals: 

\noindent $\mathrm{(O1)}$ Instead of waiting for the underlying systems to converge to intermediate reference levels $\{\bu_i[t_k]\}_{i \in \cN_D}$, steps~\eqref{eq:dual_ascent_sys}--\eqref{eq:primal_V_sys} are performed continuously (within the limits of affordable computational burden); i.e., at each instant $t_k$, one may have that $\lim_{t \rightarrow t_{k}^-} \|\by_i(t) - \bu_i[t_k]\| \neq 0$ for at least one dynamical system. This scenario is particularly relevant since step~\eqref{eq:primal_y_sys} is computationally light: it affords a closed-form solution when the inverter is operated under $\mathrm{(c1)}$ and $\mathrm{(c2)}$, and it involves a projection onto the inverter operating region under $\mathrm{(c3)}$ [cf. Figure~\ref{Fig:OIDregions}].

\noindent $\mathrm{(O2)}$ The computational time required to solve the SDP problem~\eqref{eq:primal_V_sys} is typical higher than that required by the projection operation~\eqref{eq:primal_y_sys}; especially when~\eqref{eq:primal_y_sys} affords a closed-form solution (see e.g.,~\cite{Nesterov94}, and pertinent references therein). Thus, convergence of the system outputs is investigated for the case where the update of the input reference levels $\{\bu_i[t_k]\}_{i \in \cN_D}$ and the dual variables $\{\blambda_i[t_k]\}_{i \in \cN_D}$ is performed at a faster rate than~\eqref{eq:primal_V_sys}. 

To this end, suppose that the computational time required to update matrix $\bV$ spans $M < + \infty$ time intervals; that is, if the computation of~\eqref{eq:primal_V_sys} starts at time $t_k$ based on the most up-to-date dual variables $\{\blambda_i[t_k]\}_{i \in \cN_D}$, its solution becomes available only at time $t_{k+M}$. In contrast, the controller affords the computation of steps~\eqref{eq:primal_y_sys} and~\eqref{eq:dual_ascent_sys} at each time $\{t_{k}\}_{k \in \mathbb{N}}$. To capture this asynchronous operation, consider the mapping 
\begin{align}
\label{eq:time_updates_V}
c(k) := M \left \lfloor \frac{k}{M} \right \rfloor \hspace{.5cm} k \in \mathbb{N} . 
\end{align}
Using~\eqref{eq:time_updates_V}, steps~\eqref{eq:dualsubgradient_system} for all $i \in \mathcal{N}_D$ are modified as:
\begin{subequations}
\label{eq:dualsubgradient_system_as}
\begin{align}
&  \hspace{.15cm} \by_i[t_{k}] =  \br_i \Big( \bx_i(t_{k}), \bd_i \Big)  \label{eq:sys1-obs_sys_as}  \\
& \hspace{-.2cm}  \blambda_i[t_{k+1}]   = \blambda_i[t_{k}]  + \alpha_{k+1} \Big( \bh_i(\bV[t_{c(k)}]) - \by_i[t_{k}] + \bd_i  \Big)   \label{eq:dual_ascent_sys_as} \\
& \hspace{-.2cm}  \bu_i [t_{k+1}]  = \arg \min_{\bu_i \in \cY_i} G_i(\bu_i) - \blambda_i^\sfT[t_{k+1}] \bu_i \, \label{eq:primal_y_sys_as} 
\end{align}
for all $t_k$, $k \in \mathbb{N}$. Further, matrix $\bV[t_{c(k) }] $ is updated (at the possibly slower rate) as: 
\begin{align}
\hspace{-.2cm} \bV[t_{c(k)}]  = \arg \min_{\bV \in \cV} \, H(\bV) +  \sum_{i \in \cN_D} \blambda_i^\sfT[t_{c(k)}] \, \bh_i(\bV) \label{eq:primal_V_sys_as} . \hspace{-.2cm} 
\end{align}
\end{subequations}
Since $c(k) = k$ over the interval $\{t_k, \ldots, t_{k+M-1}\}$,~\eqref{eq:primal_V_sys_as} indicates that $\bV$ is being updated  every $M$ time slots. The block diagram for~\eqref{eq:dualsubgradient_system_as} can be readily obtained by replacing step~\eqref{eq:primal_V_sys} with~\eqref{eq:primal_V_sys_as}, as well as~\eqref{eq:dual_ascent_sys} and~\eqref{eq:primal_y_sys}  with ~\eqref{eq:dual_ascent_sys_as} and~\eqref{eq:primal_y_sys_as}, respectively, in Figure~\ref{Fig:diagram}.

In the following, convergence of the system outputs to the solution of the steady-state optimization problem $\mathrm{(P1)}$ is established when the reference signals are produced by~\eqref{eq:dualsubgradient_system_as}. Of course, by setting $M = 1$, steps~\eqref{eq:dualsubgradient_system_as} coincide with~\eqref{eq:dualsubgradient_system}, and therefore the convergence claims for this more general setting naturally carry over to the synchronous setup in~\eqref{eq:dualsubgradient_system}.  

For brevity, collect the system outputs in the vector $\by := [\by_1^\sfT, \ldots, \by_{N_D}^\sfT]^\sfT$, and the dual variables  
in $\blambda := [\blambda_1^\sfT, \ldots, \blambda_{N_D}^\sfT]^\sfT$. 
In the following, it will be shown that~\eqref{eq:dual_ascent_sys} and~\eqref{eq:dual_ascent_sys_as} are in fact $\epsilon$-subgradient steps~\cite[Proposition~2]{Larsson03}  whenever $\lim_{t \rightarrow t_{k}^-} \|\by(t) - \bu[t_k]\| \neq 0$ and/or $M > 1$. Before elaborating further on the error $\epsilon[t_k]$, notice that from the compactness of sets $\cV$ and $\{\cY_i\}_{i \in \cN_D}$, it follows that there exists a constant $0 \leq G \leq + \infty$ such that the following holds:  
\begin{align}
\| \bh(\bV) - \by + \bd \|_2 \leq G \, , \quad \forall \,\, \bV \in \cV, \, \forall \,\, \by \in \cY   \label{eq:bounded_gradient_2}
\end{align} 
with $\cY := \cY_1 \times \cY_2 \times \ldots \times \cY_N$. Furthermore, given the Lipschitz-continuity of the contraction mapping~\eqref{eq:primal_y_sys}~\cite{Hong13}  
\begin{align}
\label{eq:function_u}
\bu_i(\blambda) = \arg \min_{\bu_i \in \cY_i} G_i(\bu_i) - \blambda^\sfT \bu_i , \quad  \forall i \in \cN_D 
\end{align} 
there exists $\tilde{\blambda}[t_{k}]$ satisfying
\begin{align}
\by_i[t_{k}] = \arg \min_{\bu_i \in \cY_i} G_i(\bu_i) - \tilde{\blambda}_i^\sfT[t_{k}] \bu_i , \quad  \forall i \in \cN_D 
\label{eq:lambdatilde}
\end{align} 
that is, $\by_i[t_{k}]$ would be obtained by minimizing the Lagrangian $L(\bV, \bu, \tilde{\blambda}[t_{k}])$ when $\tilde{\blambda}[t_{k}] := [\tilde{\blambda}_1^\sfT[t_{k}], \ldots, \tilde{\blambda}_{N_D}^\sfT[t_{k}]]^\sfT$ replaces $\blambda[t_{k}]$.
The following will be assumed for $\tilde{\blambda}[t_{k}]$.
\vspace{.1cm}

\begin{assumption}
\label{ass:bounds_error} 
There exists a scalar $\tilde{G}$, $0 \leq \tilde{G} < + \infty$, such that the following bound holds for all $t_k$, $k \geq 1$
\begin{align}
\| \blambda[t_k] - \tilde{\blambda}[t_k]\|_2 \leq \tilde{G}  \|\blambda[t_k] - \blambda[t_{k-1}]\|_2 . \label{eq:bounded_gradient_error}
\end{align} 
\end{assumption}
\vspace{.1cm}
Condition~\eqref{eq:bounded_gradient_error} implicitly bounds the reference signal tracking error $\|\by[t_k] - \bu[t_k]\|_2$, as specified in the following lemma. 

\vspace{.2cm}

\begin{lemma}
\label{lemma:bounded-error}
Under Assumption~\ref{ass:bounds_error}, it follows that the tracking error $\|\by[t_k] - \bu[t_k]\|_2$, $k \in \mathbb{N}$, can be bounded as 
\begin{align}
\|\by[t_k] - \bu[t_k]\|_2 \leq L \tilde{G} G \alpha_k \label{eq:bounded_tracking_error}
\end{align}
where $L$ is the Lipschitz constant of function $\bu_i(\blambda)$ in~\eqref{eq:function_u}.
\emph{Proof.} See Appendix~\ref{sec:proof-lemma-bounded-error}. $\hfill \Box$
\end{lemma}

\vspace{.1cm}

It can be noticed from~\eqref{eq:bounded_tracking_error} that the tracking error is allowed to be arbitrarily large, but 
the outputs $\by[t_k]$ should eventually follow the reference signal $\bu[t_k]$. In fact, since the sequence $\{\alpha_k\}$ is majorized by $\{\eta_k\}$, and $\eta_k \downarrow 0$, it follows that $\|\by[t_k] - \bu[t_k]\|_2 \rightarrow 0$ as $k \rightarrow \infty$. Based on this assumption, two results that establish convergence of the overall system are in order: \emph{Lemma~\ref{lemma:epsilon}} provides an analytical characterization of the $\epsilon$-subgradient step, while \emph{Theorem~\ref{thm:convergence_as}} establishes asymptotic convergence of the output powers to the optimal solution of~$\mathrm{(P1)}$.

\vspace{.1cm}

\begin{lemma}
\label{lemma:epsilon}
Suppose that at least one of the following statements is true: i) $M > 1$; ii)  at time $t_k$, $\by_i[t_{k}] \neq \bu_i[t_{k}]$ for at least one dynamical system. Then, $\bh(\bV[t_{c(k)}]) - \by[t_k] + \bd)$ is an $\epsilon$-subgradient of the dual function at $\blambda[t_k]$.  In particular, under \emph{Assumption~\ref{ass:bounds_error}} and with $M <  + \infty$, it holds that 
\begin{subequations}
\label{eq:eps_subgradient_all}
\begin{align}
& \hspace{-.3cm} \left(\bh(\bV[t_{c(k)}]) - \by[t_k] +  \bd) \right)^\sfT (\blambda - \blambda[t_k])   \nonumber \\
& \hspace{2.5cm} \geq q(\blambda) - q(\blambda[t_k]) - \epsilon[t_k] \,\,\,\, \forall \,\, \blambda \,  \label{eq:eps_subgradient}
\end{align} 
where the error $\epsilon[t_k] \geq 0$ can be bounded as 
\begin{align}
\epsilon[t_k] \leq 2 \alpha_k \tilde{G} G^2 + 2 G^2 \sum_{h = 1}^{k - c(k)} \alpha_{k-h+1} \, .\label{eq:eps_subgradient_as_2}
\end{align}
\end{subequations}
\emph{Proof.}  See Appendix~\ref{sec:proof-lemma-epsilon}. \hfill $\Box$
\end{lemma}

\vspace{.2cm}

\begin{theorem}
\label{thm:convergence_as}
Under \emph{Assumptions~\ref{ass:DynSystems}--\ref{ass:bounds_error}}, and for any $1\leq M < + \infty$,  the following holds for the closed-loop system~\eqref{eq:dualsubgradient_system_as} when a stepsize sequence $\{\alpha_k\}_{k \in \mathbb{N}}$ satisfying conditions $\mathrm{(s1)}$-$\mathrm{(s3)}$ is utilized: 

\noindent (i) $\blambda_i[t_k] \rightarrow \blambda_i^{\mathrm{opt}}$ as $k \rightarrow \infty, \forall i \in \cN_D$;  

\noindent {(ii)} $\bV[t_{c(k)}] \rightarrow \bV^{\mathrm{opt}}$  and $\{\bu_i[t_k] \rightarrow \bu_i^{\mathrm{opt}}\}_{i \in \cN_D}$ as $k \rightarrow \infty$;  

\noindent {(iii)} $\by_i(t) \rightarrow \bu_i^{\textrm{opt}}$ as $t \rightarrow \infty$, $\forall i \in \cN_D$. 

\noindent Statements (i)--(iii) hold for any initial conditions $\bV[0], \{\bu_i[0]\}_{i \in \cN_D}, \{\by_i(0)\}_{i \in \cN_D}, \{\blambda_i [0]\}_{i \in \cN_D}$, and any duration of the intervals $0 < t_{k} - t_{k-1} < \infty$, $k \in \mathbb{N}$.

\emph{Proof.} See Appendix~\ref{sec:proof-thm-convergence_as}. \hfill $\Box$
\end{theorem}

\vspace{.2cm}
 
\noindent \emph{Remark (convex relaxation or approximation of the OPF)}.  
For illustration purposes, the SDP relaxation technique for the OPF task is considered in this paper. However, the synthesis procedure outlined in the next section to develop feedback controllers that drive the inverter outputs  to solutions of pertinent convex optimization problems can be utilized in a variety of different setups. For example, it can be utilized when second-order cone relaxations~\cite{Farivar13} or linear approximations~\cite{Coffrin14,Bolognani15,DG-twrpd14} of the OPF problem are utilized. The paper considers convex relaxations or approximations of the OPF problem because dual $\epsilon$-subgradient-type methods are guaranteed to converge to optimal dual and primal solutions when applied to convex problems. The design of feedback controllers  in the case of \emph{non-convex} programs and their convergence will be the subject of future efforts.
\vspace{.1cm}

\noindent \emph{Remark (discrete variables)}. The OPF formulation considered in this paper does not include the optimization of the transformer taps at the substation as well as taps of  capacitor banks. Rather, these quantities are considered as inputs to the OPF problem, and are utilized to set the voltage at the substation~\cite{Khodr07, Paudyal11, Farivar12,OID} and form the (time-varying) admittance matrix $\bY$ in~\eqref{Pmg}. This strategy ensures full interoperability of the proposed controllers with legacy switchgear. However, it is worth noticing that transformer taps can be included in the optimization procedure by following the relaxation method investigated in e.g.,~\cite{Robbins15}.

\begin{algorithm}[t]
\label{alg:inverter-steps}
\caption{\textbf{Distributed  architecture: inverter operation}} 
\begin{algorithmic}
\small{
\FOR {$k = 1, 2, 3, \ldots$} 

\STATE [S1] Sample the inverter output $\by_i[t_k] = [P_{i}[t_k],Q_{i}[t_k]]^\sfT$ . \\

\STATE [S2] Receive $\bh_i(\bV[t_{c(k)}])$ from utility, if available (i.e., if $c(k) = k$). \\

\STATE [S3] Compute stepsize $\alpha_{k+1}$, and update $\blambda_i[t_{k+1}]$ via~\eqref{eq:dual_ascent_sys_as}.

\STATE [S4] Update the setpoints $\bu_i[t_{k+1}]$ via~\eqref{eq:primal_y_sys_as}, and implement $\bu_i[t_{k+1}]$ at the inverter.  \\

\STATE [S5] Transmit $\blambda_i[t_{k+1}]$ to the utility. \\ 

\STATE Go to step 1.

\ENDFOR

}
\end{algorithmic}
\end{algorithm}

\begin{algorithm}[t]
\label{alg:utility-steps}
\caption{\textbf{Distributed  architecture: DSO operation}} 
\begin{algorithmic}
\small{

\FOR {$k = M, 2M, 3M, \ldots$} 

\STATE [S1] Transmit $\bh_i(\bV[t_{c(k)}])$ to inverter $i$. Repeat for all  $i \in \cN_D$. \\

\STATE [S2]  Receive $\blambda_i[t_{c(k)+1}]$ from inverter $i$. Repeat for all  $i \in \cN_D$. \\
 
\STATE [S2]  Start the update of $\bV$ via~\eqref{eq:primal_V_sys_as}.

\ENDFOR

}
\end{algorithmic}
\end{algorithm}

\section{Distributed implementation of the controllers}
\label{sec:ApplicationExample}

When applied to the PV-inverter regulation problem outlined in Section~\ref{sec:Regulation}, the controller~\eqref{eq:dual_ascent_sys_as}--\eqref{eq:primal_V_sys_as} endows each PV-inverter $i \in \cN_D$ with the capability 
of steering its power output $\by_i(t) = [P_{i}(t),Q_{i}(t)]^\sfT$ towards the solution $\bu_i^{\mathrm{opt}} = [{P}_{i}^{\mathrm{opt}},{Q}_{i}^{\mathrm{opt}}]^\sfT$ of the formulated AC OPF problem. Claims (i)--(iv) of \emph{Theorem~\ref{thm:convergence_as}} hold for any duration $0 < t_{k} - t_{k-1} < \infty$, $k \in \mathbb{N}$, for any size of the distribution network.

The feedback controller~\eqref{eq:dualsubgradient_system_as} affords a \emph{distributed} implementation, where optimization tasks are distributed between the DSO and individual PV systems; see also Figure~\ref{Fig:diagram}. In particular: 

\noindent \emph{i)} Updates~\eqref{eq:dual_ascent_sys_as}-\eqref{eq:primal_y_sys_as} are implemented \emph{at each individual PV system} (they are either embedded in the inverter microcontroller, or, at the gateway level), and $\bu_i$ and $\blambda_i$ are stored locally at the same inverter; these steps are performed continuously, within affordable computational and hardware limits. Particularly,~\eqref{eq:primal_y_sys_as} is performed with the goal of pursuing inverter-related optimization objectives such as minimization of the real power curtailed~\cite{OID}. 

\noindent \emph{ii)} \emph{at the DSO}, updates~\eqref{eq:primal_V_sys_as} are performed with the goal of pursuing system-wide optimization objectives such as minimization of power losses and voltage regulation (this step is performed every $M$ time steps). 

To exchange relevant control signals, a bidirectional message passing between DSO and individual PV systems is necessary. This entails the following message exchanges every $M$ time slots: $\bh_i(\bV[t_k])$ is sent from the DSO to inverter $i$; subsequently, the up-to-date dual variable $\blambda_i[t_k]$ is sent from inverter $i$ to the DSO. 
Notice that customer $i \in \cN_D$ does \emph{not} share load demand and PV-related information with the DSO; in fact, information about the loads is not necessary when computing the update~\eqref{eq:primal_V_sys_as} at the DSO.  
Exchanging just Lagrange multipliers rather than power iterates ensures a privacy-preserving operation. The operating principles at both inverter and DSO are tabulated as Algorithm~1 and Algorithm~2, respectively, and schematically illustrated in Figure~\ref{Fig:diagram}. In Algorithm~1, is it also shown that the stepsize sequence $\{\alpha_k\}_{k \in \mathbb{N}}$ satisfying conditions $\mathrm{(s1)}$-$\mathrm{(s3)}$ is computed at the inverters' side; when changes in the load and solar irradiation conditions occur (that is, the inputs of the underlying OPF task change), inverters  exchange information to restart the sequence. For example, each inverter can utilize the sequence $\alpha_k = c/\sqrt{k-n}$, with $k \geq 1$, $c>0$ a given constant, and $n$ the index of the instant $t_n$ with the last change in the operating conditions.

Before proceeding, it is worth reiterating the underlying difference between distributed optimization approaches~\cite{Tse12,OID_TEC, Robbins15,Baldick99,Hug09,Erseghe14} and the proposed idea. In particular:

\noindent \emph{Conventional distributed optimization}: distributed OPF approaches involve the computation of steps  similar to (16), schematically illustrated in Figure~\ref{Fig:diagram}(a). Particularly, the optimal reference signals $\mathbf{u}^{\textrm{opt}}_i$ are commanded to the PV-inverters only \emph{after iterates (16) have converged}. Accordingly, the reference signals are updated at a slow time scale, dictated by the time required to solve the  distributed OPF solver.

\noindent \emph{Proposed scheme}: as illustrated in Figure~\ref{Fig:diagram}(b), the proposed controllers continuously  update the setpoints, based on current system outputs, as well as solar irradiance and load conditions. Hence, the setpoints are updated at a significantly faster rate, that may be on the same order of the inverter dynamics. As a result, the proposed controllers dynamically refresh the OPF-based targets every time that there is a variation in loads, conventional generators, and solar irradiance, and enable adaptability to fast-changing conditions.

To implement the proposed architecture, each controller at node $i \in \cN_D$ needs to collect at each time $t_k$ measurements of the demanded loads $P_{\ell, i}$ and $Q_{\ell,i}$, as well as the prevailing solar irradiation conditions (which translate to the maximum available real power). On the other hand, to perform step~\eqref{eq:primal_V_sys_as}, the DSO requires knowledge of the system topology and the load demand at nodes $i \in \cN \backslash \cN_D$; of course, any AC OPF formulation has similar prerequisites in terms of required data and measurements~\cite{OID,Farivar12,Robbins15,Bolognani15,Coffrin14,Khodr07,Paudyal11}. Since functions $\{\bh_i(\bV)\}_{i \in \cN_D}$ are linear in $\bV$, the prerequisite~\eqref{eq:cq} solely depends on the topology of the distribution network; thus,~\eqref{eq:cq} can be checked at  the utility side once matrix $\bY$ is available. 

Finally, it is worth noticing that consensus-based techniques can be adopted to speed up the computation of step~\eqref{eq:primal_V_sys_as} and improve scalability with respect to the distribution-system size~\cite{Tse12,OID_TEC,Dallanese-TSG13,Robbins15}. For example, by leveraging relevant matrix completion arguments~\cite{Dallanese-TSG13}, Lagrangian decomposition and dual gradient techniques can be adopted to decompose the computation of~\eqref{eq:primal_V_sys_as} across lines and/or portions of the system~\cite{Tse12,Robbins15}. The resultant algorithm would be similar to the one in~\eqref{eq:dualsubgradient_system_as}, but with step~\eqref{eq:primal_V_sys_as} replaced by multiple sub-problems (one per line or portion of the system) that are solved in parallel, followed by relevant dual updates; see e.g.,~\cite{Tse12,Robbins15} and~\cite{OID_TEC,Dallanese-TSG13}. It can be readily shown that the convergence claims of \emph{Theorem~\ref{thm:convergence_as}} carry over to this setup.

\section{Test cases}
\label{sec:NumericalResults}

The proposed PV-inverter control scheme is tested using a modified version of the IEEE 37-node test feeder and the IEEE 123-node test feeder.  The modified network is obtained by considering a single-phase equivalent; line impedances, shunt admittances, as well as active and reactive loads are adopted from the respective dataset.\footnote{Available at: \texttt{ewh.ieee.org/soc/pes/dsacom/testfeeders}.} The solver SDP3 is utilized to solve relevant SDPs in \texttt{MATLAB}, whereas the update of the inverter setpoints is computed in closed form. The objective of the test cases is to numerically corroborate the claims (i)--(iii) of \emph{Theorem~\ref{thm:convergence_as}}.
 
\begin{figure}[t]
\begin{center}
\includegraphics[width=7.0cm ]{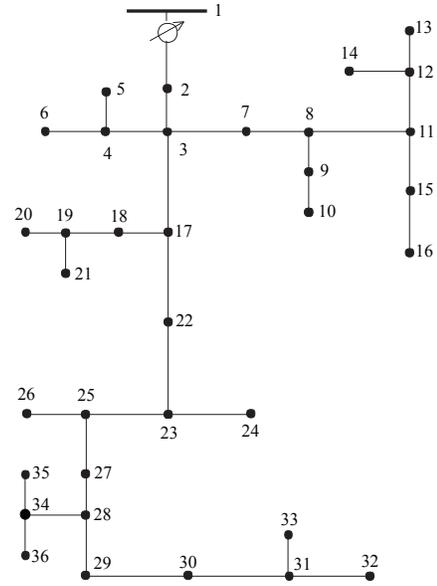}
\end{center}
\vspace{-.3cm}
\caption{IEEE 37-node test feeder considered in the test cases.}
\label{Fig:feeder}
\vspace{-.3cm}
\end{figure}

In the OPF problem, the voltage limits $V^\textrm{min}$ and $V^\textrm{max}$ are set to $0.95 \mathrm{pu}$ and $1.05 \mathrm{pu}$, respectively. In the first test, the IEEE 37-node test feeder illustrated in Figure~\ref{Fig:feeder} is utilized. The voltage magnitude at the point of common coupling is fixed to $|V_0| = 1$ pu, and it is  presumed that $6$ PV systems are present in the network and they are located at nodes $4, 11, 22, 26, 29$ and $32$. Following the technical approach of~\cite[Ch.~8]{Iravanibook10} and~\cite{Irminger12}, a first-order system is utilized to model the real and reactive power dynamics of each inverter. Further, inverters implement strategy $\mathrm{(c3)}$, and their regions of possible operating points is formed based on the inverter power ratings $\{S_i\}_{i \in \cN_D}$ and the available active powers $\{P_i^{\textrm{av}}\}_{i \in \cN_D}$. Specifically, the power ratings are assumed to be $50, 120, 50, 100, 120,$ and $80$ kVA, whereas the following values for the available powers $\bp^{\mathrm{av}} := [P_1^{\textrm{av}}, \ldots, P_{N_D}^{\textrm{av}}]^\sfT$ are considered in order to test the adaptability of the feedback controller to changing prevailing conditions (with time intervals normalized with respect to the time constant $\tau$): 

\noindent $\mathrm{(I1)}$ $\bp^{\mathrm{av}}(t) = [22, 67, 21, 50, 68, 40]^\sfT$ $\mathrm{kW}$, $t/\tau \in [1,200]$; 

\noindent $\mathrm{(I2)}$ $\bp^{\mathrm{av}}(t) = [25, 80, 24, 55, 85, 45]^\sfT$ $\mathrm{kW}$, $t/\tau \in [201, 400]$;  

\noindent $\mathrm{(I3)}$ $\bp^{\mathrm{av}}(t) = [31, 92, 29, 65, 92, 54]^\sfT$ $\mathrm{kW}$, $t/\tau \in [401,600]$;  

\noindent $\mathrm{(I4)}$ $\bp^{\mathrm{av}}(t) = [26, 84, 25, 57, 86, 47]^\sfT$ $\mathrm{kW}$, $t/\tau \in [601,700]$. 

\noindent At $t = 0$, the output active and reactive powers are $0$ kW and $0$ kVAr, respectively. No minimum power factor constraints are enforced (i.e., is $\theta = \pi/2$), and ${P}_{i}^{\textrm{min}}$ is set to $0$~\cite{Tonkoski11}. In this fiesr test,  $H(\bV)$ models the cost of power drawn from the substation as $H(\bV) =  (\trace(\bPhi_0 \bV))^2 + 10 \times \trace(\bPhi_0 \bV)$. On the other hand, the function $G_i(P_i, Q_i)$ is set to 
\begin{align}
G_i(P_i, Q_i) & =  a_i (P_{i}^{\textrm{av}} - P_i)^2 + b_i (P_{i}^{\textrm{av}} - P_i) \nonumber \\
& \hspace{.5cm} + c_i Q_i^2 + d_i |Q_i|
 \label{eq:Gi_test}
\end{align} 
in order to minimize the amount of curtailed real power, as well as the amount of reactive power provided. It is, however, worth emphasizing that various alternative cost functions can be accommodated in the proposed framework, and~\eqref{eq:Gi_test} is utilized as a representative example. Coefficients $a_i, b_i, c_i, d_i$ are set to $a_i = 1, b_i = 10, c_i = 0.01, d_i = 0.01$ for  $i = 1,\ldots, 4$, and $a_i = 1, b_i = 10, c_i = 0.03, d_i = 0.03$ for  $i = 5, 6$.  With this setup, the SDP relaxation was first tested with these input data, the SDP solver identified solutions with rank-$1$ matrices $\bV^{\mathrm{opt}}$~\cite{Bai08,LavaeiLow}.

At each inverter $i \in \cN_D$, the reference signal $\bu_i[t_k]$ is updated every $t = \tau$ $\mathrm{sec}$; i.e., $t_{k} - t_{k-1}  = \tau$ for all $k \in \mathbb{N}$. This implies that a new reference signal $\bu_i[t_k]$ is applied to each inverter  faster than the output power settling time (which corresponds to approximately $5 \tau$ for a first-order system). On the other hand, matrix $\bV[t_k]$ is updated every $t = 2 \tau\,\mathrm{sec}$; i.e., $M = 2$ in~\eqref{eq:dualsubgradient_system_as}. This means that the inverter setpoints $\{\bu_i[t_k]\}_{i \in \cN_D}$ are updated at a faster rate than matrix $\bV[t_k]$. The stepsize in~\eqref{eq:dualsubgradient_system_as} is set to $\alpha_k = 4/\sqrt{k-n}$, with $k \geq 1$ and $n$ the index of the instant $t_n$ with the last step change.

\begin{figure}[t]
\begin{center}
\subfigure[]{\hspace{-.3cm}\includegraphics[width=9.3cm]{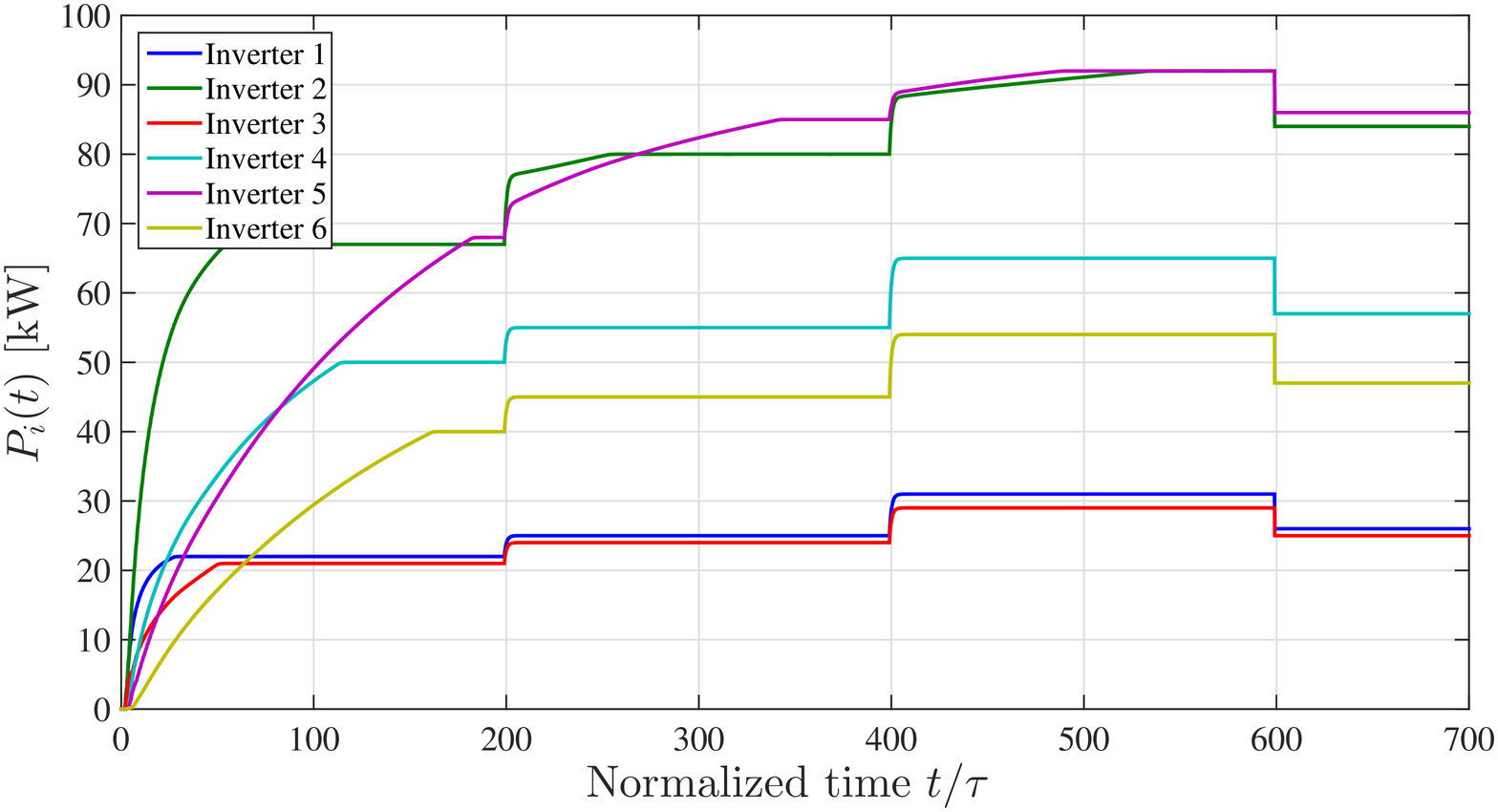}} 
\subfigure[]{\hspace{-.3cm}\includegraphics[width=9.3cm]{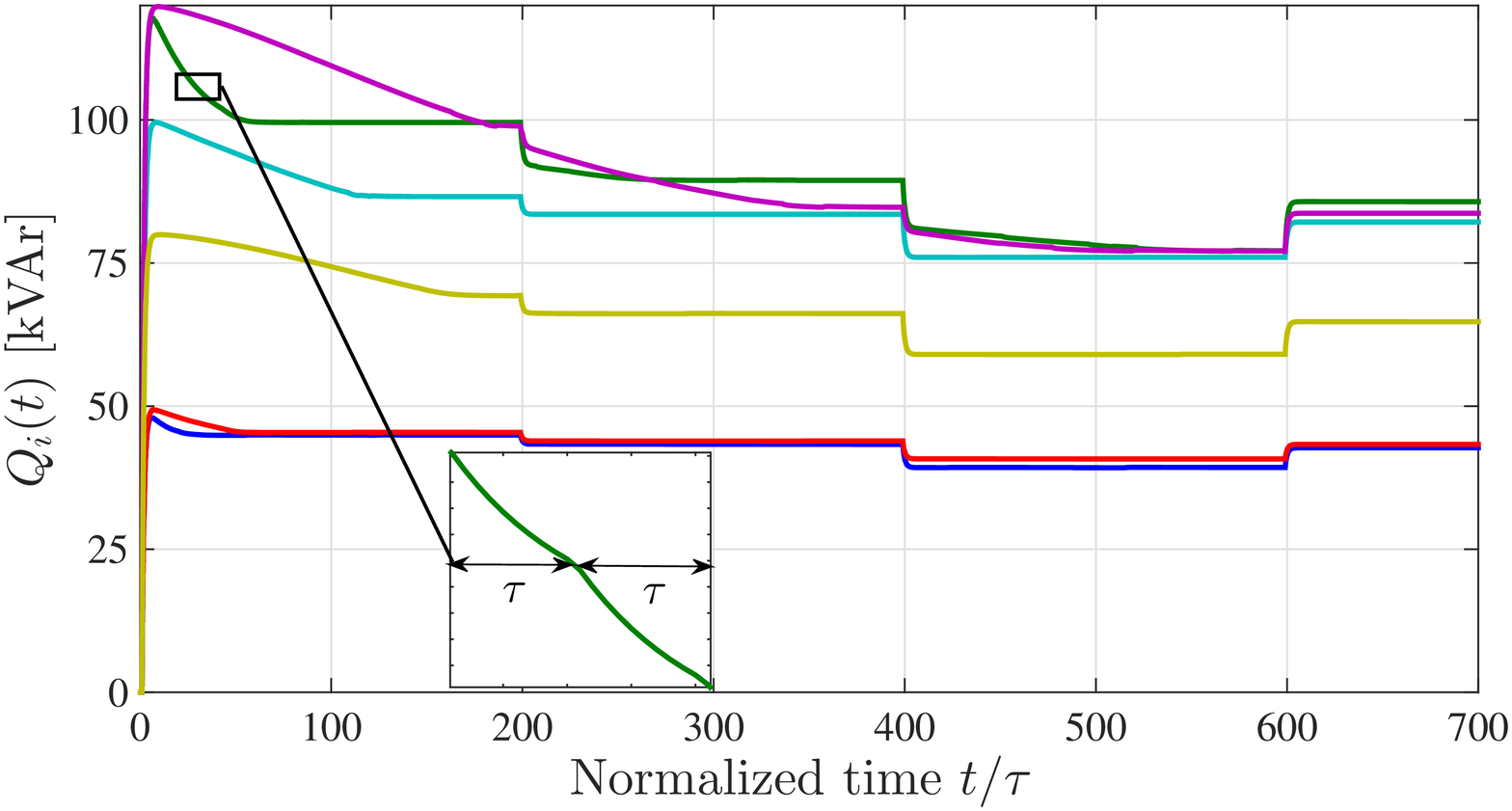}}  
\end{center}
\vspace{-.3cm}
\caption{Convergence of~\eqref{eq:dualsubgradient_system_as}, when the inverter-power dynamics are approximated as first-order systems with time constant $\tau$, for four different solar irradiance conditions. Plots illustrate the convergence of the real and reactive powers to the solutions of the formulated OPF problem. }
\label{Fig:convergence}
\vspace{-.5cm}
\end{figure}

Figure~\ref{Fig:convergence} illustrates the evolution of the real and reactive powers generated by the inverters. It can be seen that the inverter outputs  $\{\by_i[t_k] = [P_i(t), Q_i(t)]^\sfT\}$ converge in all the considered intervals $\mathrm{(I1)}$--$\mathrm{(I4)}$, and the output powers at convergence coincide with the solutions of the OPF $\mathrm{(P1)}$; for example, before the step change at $t = 200 \tau$, the active and reactive powers converged to the OPF solution $21.8, 66.9, 20.9, 67.9, 39.9$ kW and $39.2, 85.6, 40.7, 77.1, 31.4, 39.8$ kVAr. This corroborates the claims of Theorem~\ref{thm:convergence_as}. Figure~\ref{Fig:convergence}(b) also provides a snapshot of the evolution of the output reactive power for inverter $2$; it can be seen that a new reference level is applied after $\tau$ seconds, before $Q_2(t)$ settles around the intermediate setpoint. It is also interesting to note that, in the considered setup, the steady-state reactive powers coincide with the available powers $\bp^{\mathrm{av}}(t)$, and reactive compensation turns out to be the optimal ancillary service strategy. Similar trajectories would have been obtained when the loads are also time varying. Future efforts will explore variations of load and solar irradiance that may have the same temporal scale of the dynamics of~\eqref{eq:dualsubgradient_system_as}. 

In the second test case, a scenario with high PV-penetration is considered. Specifically, $17$ PV systems are assumed located at nodes $4, 11, 13, 16, 17, 20, 22, 23, 26, 28, 29, 30, 31, 32, 33, 34, 36$ of the feeder depicted in Figure~\ref{Fig:feeder}, and their  AC power ratings are assumed to be $\mathbf{s} = [50, $ $200, 220, 120, 200, 120, 150, 50, 280, 100, 250, 100, 120, 200,$ $110, 250, 150]$ kVA.   Similar to the previous case, step changes in the solar irradiation (and, hence, in the available powers $\bp^{\mathrm{av}}$) are considered in order to test the adaptability of the feedback controller to changing prevailing conditions. Specifically, the following values are tested: \\
\noindent $\mathrm{(I1)}$ $\bp^{\mathrm{av}}(t) = 0.7  \, \mathbf{s}$ $\mathrm{kW}$, $t/\tau \in [1,200]$; \\
\noindent $\mathrm{(I2)}$ $\bp^{\mathrm{av}}(t) = 0.8  \, \, \mathbf{s}$ $\mathrm{kW}$, $t/\tau \in [201, 400]$;  \\
\noindent $\mathrm{(I3)}$ $\bp^{\mathrm{av}}(t) = \mathbf{s}$ $\mathrm{kW}$, $t/\tau \in [401,500]$;  \\
\noindent $\mathrm{(I4)}$ $\bp^{\mathrm{av}}(t) = 0.6 \, \mathbf{s}$ $\mathrm{kW}$, $t/\tau \in [501,700]$. \\
At each inverter $i \in \cN_D$, the reference signal $\bu_i[t_k]$ is updated every $t = \tau$ $\mathrm{sec}$, while $\bV[t_k]$ is updated every $t = 2 \tau\,\mathrm{sec}$. The voltage magnitude at the substation is set to $1.02$ pu, while coefficients $a_i, b_i, c_i, d_i$  in~\eqref{eq:Gi_test} are set to $a_i = 1, b_i = 1, c_i = 0.5, d_i = 3$ for  all $i = 1,\ldots, 17$. All the other simulation parameters are similar to the previous test case.  The SDP relaxation was tested under this setup,  and the SDP solver identified solutions with rank-$1$ matrices $\bV^{\mathrm{opt}}$ for all the four cases considered~\cite{Bai08,LavaeiLow}.

\begin{figure*}[t]
\begin{center}
\subfigure[]{\hspace{-.3cm}\includegraphics[width=9.3cm]{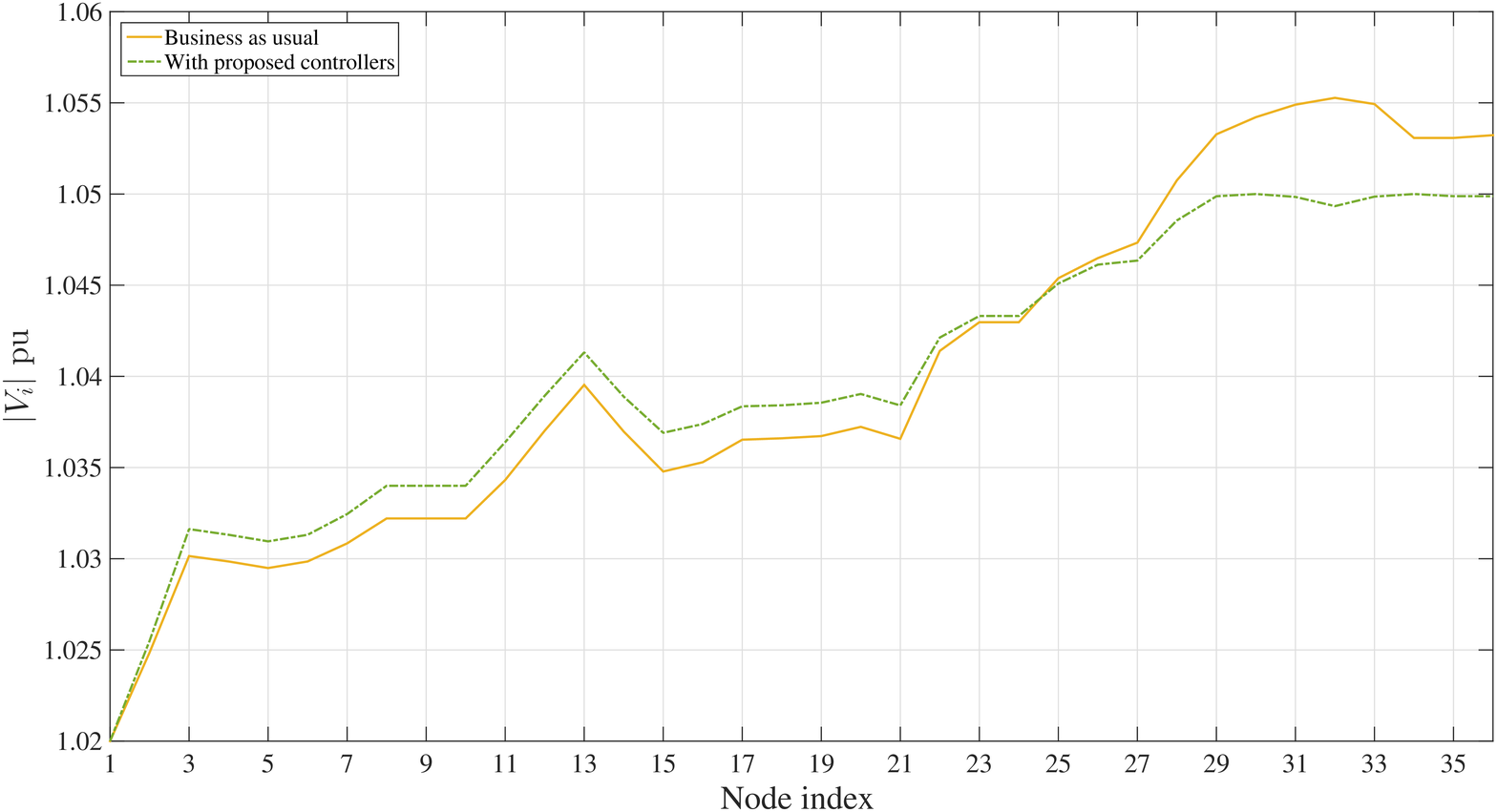}} 
\subfigure[]{\hspace{-.3cm}\includegraphics[width=9.3cm]{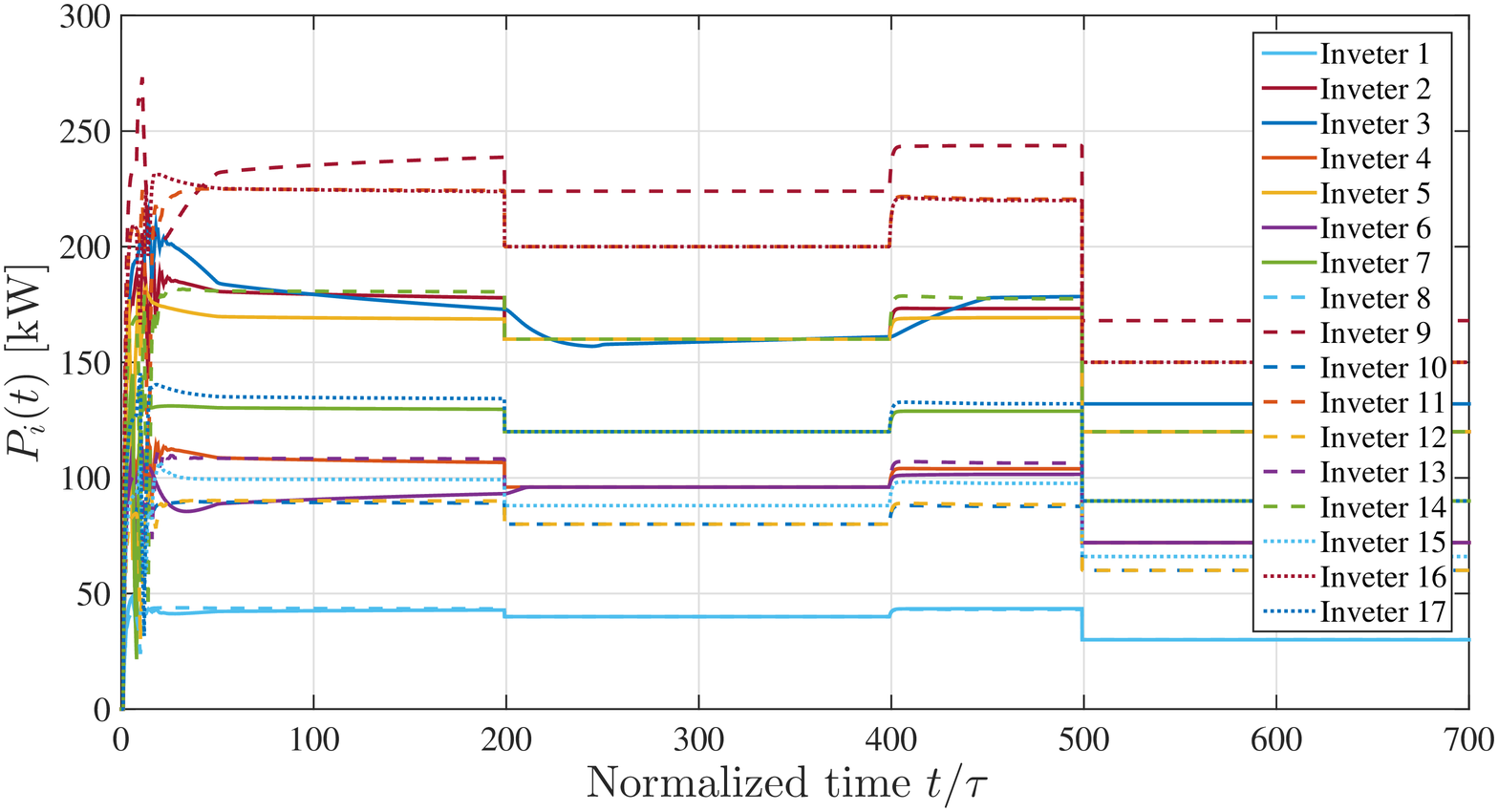}} 
\subfigure[]{\hspace{-.3cm}\includegraphics[width=9.3cm]{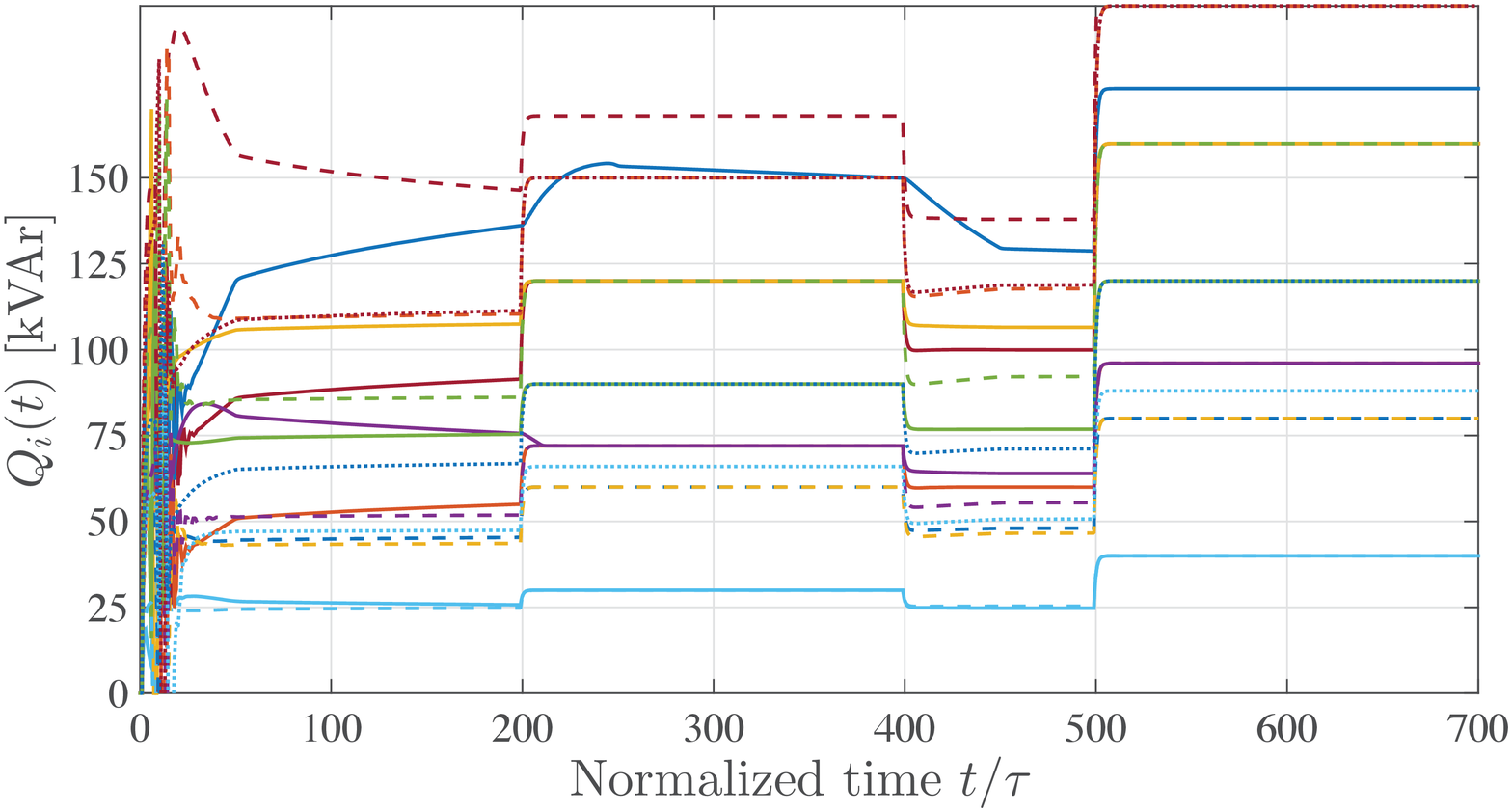}} 
\end{center}
\vspace{-.3cm}
\caption{Test case with high PV-penetration. (a) Voltage profile at $t = 450 \tau$. (b) Real powers provided by the inverters. (c) Reactive powers provided by the inverters. }
\label{Fig:convergence_2}
\vspace{-.5cm}
\end{figure*}

With this setup, when inverters operate at unity power factor and set $P_i = P^{\mathrm{av}}$, the voltage magnitudes exceed the upper limit of $1.05$ pu  during the interval  $t/\tau \in [401,500]$  in 10 nodes. Specifically, the voltage profile is shown with the yellow trajectory in Figure~\ref{Fig:convergence_2}(a).

The objective of this test case is twofold: i) demonstrate voltage regulation capabilities of the proposed scheme; and, ii) demonstrate that the convergence speed is not deteriorated when a higher number of PV systems are controlled. As for objective i), it can be clearly seen in  Figures~\ref{Fig:convergence_2}(a) that the voltages are steadily kept within the limits $V^{\textrm{min}}$ and $V^{\textrm{max}}$; particularly, the green trajectory in Figure~\ref{Fig:convergence_2}(a) shows that the proposed scheme favors voltage regulation even during peak generation conditions, while minimizing the amount of curtailed real power [cf.~\eqref{eq:Gi_test}]. Figures~\ref{Fig:convergence_2}(c) and~\ref{Fig:convergence_2}(d) illustrate  the evolution of the real and reactive powers generated by the inverters. Comparing Figure~\ref{Fig:convergence} with  Figures~\ref{Fig:convergence_2}(c) and~\ref{Fig:convergence_2}(d), it can be noticed that the proposed controllers still provide fast adaptation capabilities to changes in the solar irradiation; and furthermore, convergence speed is not degraded when an increased number of PV systems are controlled.

\begin{figure}[t]
\begin{center}
\includegraphics[width=9.0cm ]{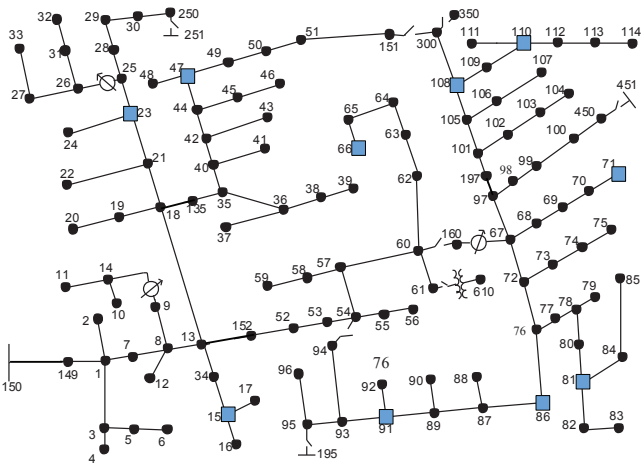}
\end{center}
\vspace{-.3cm}
\caption{IEEE 123-node test feeder considered in the third test case. Blue squares represent nodes at which PV systems are installed.}
\label{Fig:feeder123}
\vspace{-.3cm}
\end{figure}

In the third test case,  the IEEE 123-node test feeder illustrated in Figure~\ref{Fig:feeder123} is utilized. The voltage magnitude at the point of common coupling is fixed to $|V_0| = 1$ pu, and it is  presumed that $10$ PV systems are located at nodes $15, 23, 47, 66, 71, 81, 86, 91, 108$ and $110$. Inverters implement strategy $\mathrm{(c3)}$, and their AC power ratings amount to  $\mathbf{s} = [500, 450, 200, 300, 200, 200, 150, 150, 200, 350]$ kVA.  Changes in the solar irradiation are considered in order to test the adaptability of the feedback controller  for this larger system; the following values are tested: \\
\noindent $\mathrm{(I1)}$ $\bp^{\mathrm{av}}(t) = 0.8  \, \, \mathbf{s}$ $\mathrm{kW}$, $t/\tau \in [1, 200]$;  \\
\noindent $\mathrm{(I2)}$ $\bp^{\mathrm{av}}(t) = \mathbf{s}$ $\mathrm{kW}$, $t/\tau \in [201,400]$;  \\
\noindent $\mathrm{(I3)}$ $\bp^{\mathrm{av}}(t) = 0.6 \, \mathbf{s}$ $\mathrm{kW}$, $t/\tau \in [401,500]$. \\
Similar to the previous test cases, the reference signals $\{\bu_i[t_k]\}$ are updated every $t = \tau$ $\mathrm{sec}$, while $\bV[t_k]$ is updated every $t = 2 \tau\,\mathrm{sec}$. In the OPF problem, function $H(\bV)$ captures (the cost of) power losses, and it is set to $H(\bV) =  (\trace(\bL \bV))^2 + 5 \times \trace(\bL \bV)$, with matrix $\bL$ formed as described in~\cite{Dallanese-TSG13}. The coefficients in~\eqref{eq:Gi_test} are set to $a_i = 1, b_i = 10, c_i = 0.5, d_i = 3$ for  all inverters. All the other simulation parameters are similar to the previous test case.  The SDP relaxation was tested,  and the SDP solver identified solutions with rank-$1$ matrices $\bV^{\mathrm{opt}}$.

\begin{figure}[t]
\begin{center}
\subfigure[]{\hspace{-.3cm}\includegraphics[width=9.3cm]{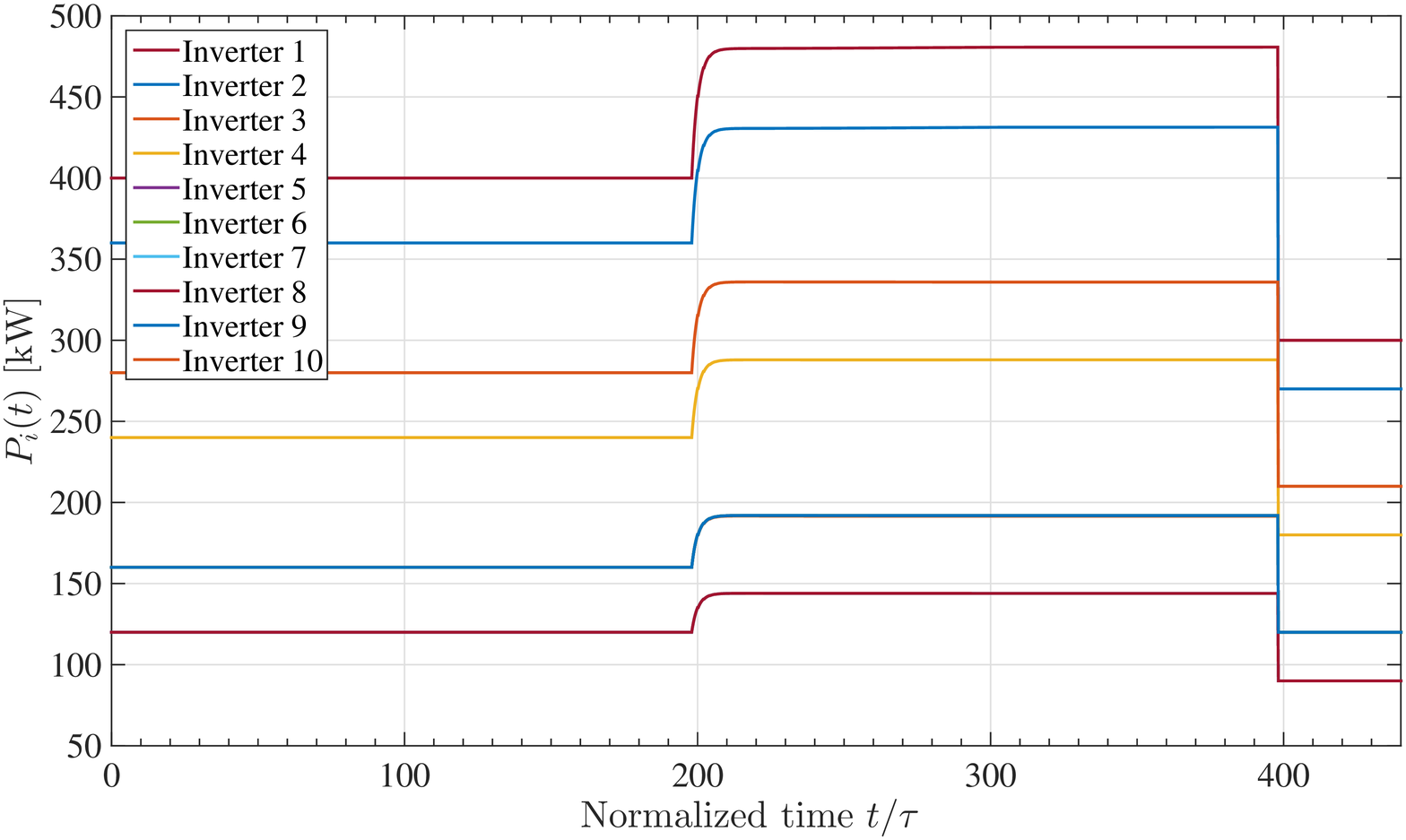}} 
\subfigure[]{\hspace{-.3cm}\includegraphics[width=9.3cm]{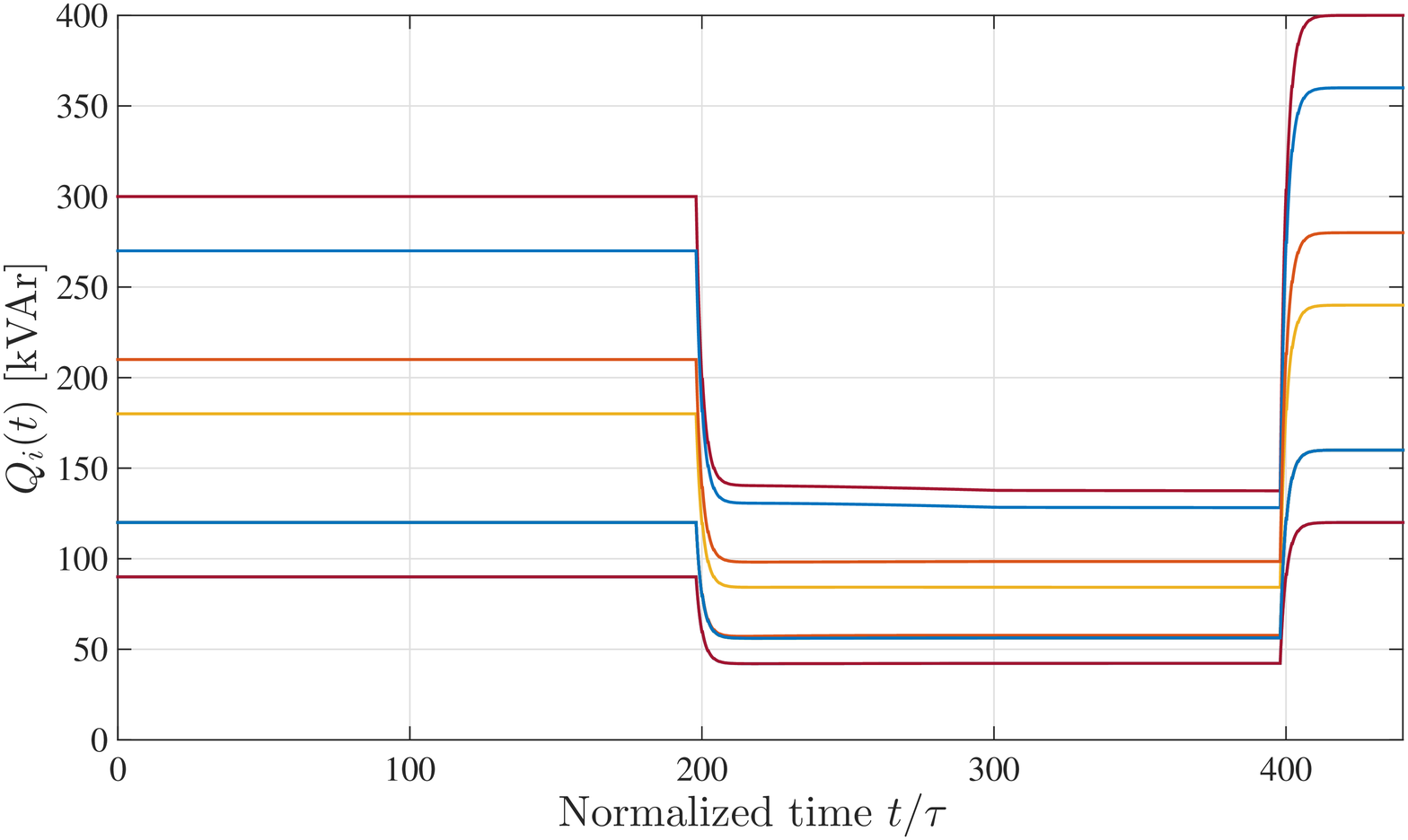}} 
\end{center}
\vspace{-.3cm}
\caption{Convergence of~\eqref{eq:dualsubgradient_system_as}, for the modified IEEE 123-node test feeder. The inverter-power dynamics are approximated as first-order systems with time constant $\tau$, for four different solar irradiance conditions. Plots illustrate the convergence of the real and reactive powers to the solutions of the formulated OPF problem. }
\label{Fig:convergence_123}
\vspace{-.5cm}
\end{figure}

Figures~\ref{Fig:convergence_123}(c) and~\ref{Fig:convergence_123}(d) illustrate  the evolution of the real and reactive powers generated by the 10 inverters. It can be seen that the inverters quickly regulate the power outputs to new OPF setpoints. In particular, comparing Figure~\ref{Fig:convergence},  Figures~\ref{Fig:convergence_2} and~\ref{Fig:convergence_123}, it can be noticed that the convergence speed of the proposed controllers is not degraded when a larger distribution network is controlled. Notice that inverters are required to curtail real power in order to adhere to voltage limits. 

Finally, it is worth mentioning that, for the update~\eqref{eq:primal_V_sys_as}, well-established complexity bounds for convex programs such as SDPs~\cite{Nesterov94} exist; these bounds quantify how the worst-case computational complexity increases with the number of variables (i.e., the network size). Further, in case of an SDP, provably convergent parallelization techniques can also be leveraged to speed up the computation of~\eqref{eq:primal_V_sys_as}; see, e.g.,~\cite{Summers_IFA_5019}.   

\section{Concluding Remarks and Future Work}
\label{sec:Conclusions}

This paper considered a distribution network featuring PV systems, and addressed the synthesis of feedback controllers that seek inverter setpoints corresponding to AC OPF solutions. To this end, dual $\epsilon$-subgradient methods and SDP relaxations were leveraged. Global convergence of PV-inverter output powers was analytically established and numerically corroborated. Although the focus was on PV systems, the framework naturally accommodates different types of inverter-interfaced energy resources. The development of provably convergenct feedback controllers that seek the solutions of non-convex OPF formulations will be the subject of future research efforts.

\appendix

\subsection{Proof of Lemma~\ref{lemma:bounded-error}}
\label{sec:proof-lemma-bounded-error}

Recall first that, given the strong convexity of $G_i(\bu_i)$, it turns out that function $\bu_i(\blambda)$ in~\eqref{eq:function_u} is Lipschitz continuous (in $\lambda$), with a constant denoted here as $L$~\cite{Hong13}. Then, recalling that $\by[t_k] = \bu_i(\tilde{\lambda}[t_k])$ and $\bu[t_k] = \bu_i(\lambda[t_k])$, it follows that the left-hand side of~\eqref{eq:bounded_tracking_error} can be bounded as
\begin{subequations}
\begin{align}
\|\by[t_k] - \bu_i[t_k]\|_2 & \leq  L \|\tilde{\blambda}[t_k] - \blambda[t_k]\|_2 \\
 & \leq L \tilde{G} G \alpha_k  \label{eq:proof-track-s3} 
\end{align}
where~\eqref{eq:proof-track-s3} is obtained by using
 the following bound (which originates from \emph{Assumption 4}): 
\begin{align}
\| \blambda[t_k] - \tilde{\blambda}[t_k]\|_2 & \leq \tilde{G}  \|\blambda[t_k] - \blambda[t_{k-1}]\|_2  \label{eq:step1} \\
&\hspace{-0.6in}  \leq \tilde G \| \alpha_{k} ( \bh(\bV[t_{c(k-1)}]) + \bg(\by[t_{k-1}], \bd)) \|_2 \label{eq:step2} \\
&\hspace{-0.6in} \leq \tilde G  G \alpha_k. \label{eq:step3}
\end{align}
Note that~\eqref{eq:step2} follows from the dual update in~\eqref{eq:dual_ascent_sys_as}, and~\eqref{eq:step3} follows from~\eqref{eq:bounded_gradient_2}. 
\end{subequations}

\subsection{Proof of Lemma~\ref{lemma:epsilon}}
\label{sec:proof-lemma-epsilon}

Recall  that $\bh(\bV[t_k]) - \bu[t_k] + \bd$ is the gradient of the dual function~\eqref{eq:dualFunction} evaluated at $\blambda[t_k]$~\cite{Bertsekas_ConvexAnalysis}. Let $\bg(\bu, \bd) := \bd - \bu$ for exposition simplicity, and consider decomposing~\eqref{eq:dualFunction} as $q(\blambda) = q_V(\blambda) + q_u(\blambda)$, with
\begin{subequations}
\begin{align}
q_V(\blambda) & := \min_{\bV \in \cV} H(\bV) +  \blambda^\sfT \bh(\bV) \label{eq:dual_function_decomp_1}, \\
q_u(\blambda) & := \min_{\bu \in \cU}  G(\bu) + \blambda^\sfT \bg(\bu, \bd) \label{eq:dual_function_decomp_2}
\end{align} 
where $G(\bu) := \sum_{i \in \cN_D} G_i(\bu_i)$.
Then, it will be shown that 
\begin{align}
& \hspace{-.3cm} \bg^\sfT(\by[t_k], \bd) (\blambda - \blambda[t_k])   \geq q_u(\blambda) - q_u(\blambda[t_k]) - \epsilon_u[t_k] \label{eq:err_sub_u}\\
& \hspace{-.3cm} \bh^\sfT(\bV[t_{c(k)}]) (\blambda - \blambda[t_k])  \geq q_V(\blambda) - q_V(\blambda[t_k]) - \epsilon_V[t_k]  \label{eq:err_sub_V}
\end{align} 
with $\epsilon_u[t_k] \leq 2 \alpha_k \tilde{G} G^2$ and $\epsilon_V[t_k] \leq 2 G^2 \sum_{h = 1}^{k - c(k)} \alpha_{k-h+1}$. 

To show~\eqref{eq:err_sub_u}, consider the gradient of $q_u(\blambda)$ evaluated at $\tilde \blambda[t_k]$, which by definition leads to the inequality $\bg^\sfT(\by[t_k], \bd) (\blambda - \tilde \blambda[t_k]) \geq q_u(\blambda) - q_u(\tilde \blambda[t_k])$ for all $\blambda$; then,  add $\bg^\sfT(\by[t_k], \bd) (\tilde \blambda[t_k] - \blambda[t_k])$ on both sides to obtain 
\small
\begin{align}
& \bg^\sfT(\by[t_k], \bd) (\blambda - \blambda[t_k]) \geq q_u(\blambda) - q_u(\tilde \blambda[t_k]) \nonumber \\
& \hspace{3cm} + \bg^\sfT(\by[t_k], \bd) (\tilde  \blambda[t_k] - \blambda[t_k]) 
\end{align}
\normalsize
and add and subtract $q_u(\blambda[t_k])$ to the right-hand-side
\begin{align}\label{eq:ineq}
&\bg^\sfT(\by[t_k], \bd) (\blambda - \blambda[t_k]) \geq \,\,q_u(\blambda) - q_u(\blambda[t_k])   \\
& +q_u(\blambda[t_k]) - q_u(\tilde \blambda[t_k]) + \bg^\sfT(\by[t_k], \bd) (\tilde \blambda[t_k] - \blambda[t_k]). \nonumber 
\end{align}
In~\eqref{eq:ineq}, define $\epsilon_u[t_k]:= q_u(\tilde \blambda[t_k]) - q_u(\blambda[t_k]) + \bg^\sfT(\by[t_k], \bd)(\blambda[t_k] - \tilde \blambda[t_k])$. 
By using the definition of the gradient of the function $q_u(\blambda)$ at $\blambda[t_k]$, and applying the Cauchy-Schwartz inequality, one has that
\begin{align}
\epsilon_u[t_k] &\leq  \bg^\sfT(\bu[t_k], \bd)(\tilde \blambda[t_k] - \blambda[t_k]) \nonumber \\ 
& \hspace{.4cm} + \bg^\sfT(\by[t_k], \bd)(\blambda[t_k] - \tilde \blambda[t_k])\label{eq:edef0}  \\
&\leq  2 G \,\, \|\tilde \blambda[t_k] - \blambda[t_k]\|_2 \leq 2 \alpha_k \tilde G G^2 \label{eq:edef1}
\end{align}
where~\eqref{eq:dual_ascent_sys_as},~\eqref{eq:bounded_gradient_2}, and~\eqref{eq:bounded_gradient_error} were used to obtain~\eqref{eq:edef1} from~\eqref{eq:edef0}. Next, to show~\eqref{eq:err_sub_V}, begin with the inequality $\bh^\sfT(\bV[t_{c(k)}]) (\blambda - \blambda[t_{c(k)}]) \geq q_V(\blambda) - q_V(\blambda[t_{c(k)}])$. Adding $\bh^\sfT(\bV[t_{c(k)}]) (\blambda[t_{c(k)}] - \blambda[t_k])$ to both sides of the inequality,
\begin{align}
\bh^\sfT(\bV[t_{c(k)}]) (\blambda - \blambda[t_k]) & \geq q_V(\blambda) - q_V(\blambda[t_{c(k)}]) \nonumber \\
&  \hspace{-1.2cm}  +  \bh^\sfT(\bV[t_{c(k)}]) (\blambda[t_{c(k)}] - \blambda[t_k]) \, .
\label{eq:eps_subgradient_as_proof_1} 
\end{align} 
Adding and subtracting the sequences $\{q_V(\blambda[t_{k-h+1}])\}_{h = 1}^{k - c(k)}$ and $\{\bh^\sfT(\bV[t_{c(k)}]) (\blambda[t_{k-h+1}])\}_{h = 2}^{k - c(k)}$ to the  right-hand-side of~\eqref{eq:eps_subgradient_as_proof_1}, and suitably rearranging terms, one obtains
  \begin{align}
\hspace{-.15cm} \bh^\sfT(\bV[t_{c(k)}]) (\blambda - \blambda[t_k]) & \geq q_V(\blambda) - q_V(\blambda[t_{k}]) - \epsilon_V[t_k] \hspace{-.2cm}
\end{align} 
where $\epsilon_V[t_k]$ is defined as 
\begin{align}
\epsilon_V[t_k] & := \sum_{h = 1}^{k - c(k)} \left(q_V(\blambda[t_{k - h}]) - q_V(\blambda[t_{k - h + 1}]) \right) \nonumber \\
& - \sum_{h = 1}^{k - c(k)} \bh^\sfT(\bV[t_{c(k)}]) (\blambda[t_{k - h}] - \blambda[t_{k - h + 1}]) \, .
\label{eq:eps_subgradient_as_proof_2} 
\end{align} 
Using the definition of the gradient, the Cauchy-Schwartz inequality, and~\eqref{eq:bounded_gradient_2},~\eqref{eq:eps_subgradient_as_proof_2} can be bounded as:
\small
\begin{align}
\epsilon_V[t_k] \leq 2 G \sum_{h = 1}^{k - c(k)} \|\blambda[t_{k - h + 1}] - \blambda[t_{k - h}] \|_2 \, .
\label{eq:eps_subgradient_as_proof_3} 
\end{align} 
\normalsize
Finally, upon using~\eqref{eq:dual_ascent_sys_as} and~\eqref{eq:bounded_gradient_2},~\eqref{eq:eps_subgradient_as_proof_3} can be further bounded as $ 2 G \sum_{h = 1}^{k - c(k)} \|\blambda[t_{k - h + 1}] - \blambda[t_{k - h}] \|_2 \leq 2 G^2 \sum_{h = 1}^{k - c(k)} \alpha_{k-h+1}$. 
\end{subequations}

\subsection{Proof of Theorem~\ref{thm:convergence_as}}
\label{sec:proof-thm-convergence_as}

\emph{Claims (i)--(ii)}. Boundedness and convergence of the dual iterates can be proved by leveraging the results in~\cite[Theorem~3.4]{Kiwiel04}. In particular, it suffices to show that the following technical requirement is satisfied in the present setup:  
\begin{subequations}
\begin{align}
\sum_{k = 0}^{+ \infty} \alpha_k \epsilon[t_k] = \sum_{k = 0}^{+ \infty} \alpha_k (\epsilon_V[t_k] + \epsilon_u[t_k]) <  + \infty \, . \label{proofT1-s1}
\end{align}
From \emph{Lemma~\ref{lemma:epsilon}}, it it can be shown that
\begin{align}
\sum_{k = 0}^{+ \infty} \alpha_k \epsilon_u[t_k] & \leq \sum_{k = 0}^{+ \infty}  2 \alpha_k^2 \tilde{G} G^2 \leq 2  \tilde{G} G^2 \sum_{k = 1}^{+ \infty}  \eta_k^2 \label{proofT1-s3}
\end{align}
where the second inequality in~\eqref{proofT1-s3} follows from the fact that $\alpha_k \leq \eta_k$ for all $k$. Since $\sum_{k = 0}^{+ \infty} \eta_k^2 < + \infty$, the series $\sum_{k = 0}^{+ \infty} \alpha_k \epsilon_u[t_k]$ is finite. As for the error $\epsilon_V[t_k]$, one has that 
\small
\begin{align}
\sum_{k = 0}^{+ \infty} \alpha_k \epsilon_V[t_k]  & \leq 2 G^2 \sum_{k = 0}^{+ \infty}  \alpha_k \sum_{h = 1}^{k - c(k)} \alpha_{k-h+1} \label{proofT2-s2} \\
&  \leq 2 G^2 \sum_{k = 0}^{+ \infty}  \alpha_k \sum_{h = 1}^{M-1} \alpha_{k-h+1} \label{proofT2-s3} \\
& \leq 2 G^2 \sum_{k = 0}^{+ \infty}  \sum_{h = 1}^{M-1} \eta^2_{k-h+1} \label{proofT2-s4}
\end{align}
\end{subequations}
\normalsize
where the fact that $\max \{k - c(k)\} = M-1$ is utilized in~\eqref{proofT2-s3}, and~\eqref{proofT2-s4} follows from~\eqref{proofT2-s3} since the sequence $\{\eta_k\}_{k \in \mathbb{N}}$ majorizes $\{\alpha_k\}_{k \in \mathbb{N}}$, and it is  monotonically decreasing.  Since the series $\{\eta^2_{k}\}$ is square-summable, $\sum_{k = 0}^{+ \infty} \alpha_k \epsilon_V[t_k]$ is finite. 

\emph{Claim (iii)} From the convexity of the Lagrangian in the primal variables, it follows that optimal primal variables can be uniquely recovered as $\{ \bV^{\mathrm{opt}}, \bu^{\textrm{opt}}\} = \arg \min_{\bV \in \cV, \bu \in \cU} L\left(\bV, \bu, \blambda^{\mathrm{opt}}\right)$. 

\emph{Claim (iv)} At convergence, the reference signal is constant, with value $\bu_i^{\textrm{opt}}$. Then, $\by_i(t) \rightarrow \bu_i^{\textrm{opt}}$ as $t \rightarrow \infty$ by~\eqref{eq:sys1-equilibrium}.

\subsection{Extension to multi-phase systems}
\label{sec:three-phase}

For notation and exposition simplicity, Sections~\ref{sec:ProblemFormulation} and~\ref{sec:ApplicationExample} considered a balanced distribution network. However, the proposed framework can be extended to multi-phase systems as detailed in~\cite{Dallanese-TSG13,Gan-PSCC14,Robbins15} and briefly explained in the following.

Define as $\cP_{ij} \subseteq \{a,b,c\}$ and $\cP_{i} \subseteq \{a,b,c\}$ the sets of phases of line $(i,j) \in \cE$ and node $i \in \cN$, respectively. Hereafter, a superscript $(\cdot)^\phi$ is utilized to assign relevant electrical quantities to a specific phase. For example, $V_i^{\phi} \in \mathbb{C}$ denotes the complex line-to-ground voltage at node $i \in \cN$ and phase $\phi \in \cP_i$, whereas $I_i^{\phi} \in \mathbb{C}$ is the phasor representation of the current injected at the same node and phase; further, $P_i^\phi$ and $Q_i^\phi$ denote the output real and reactive powers of a PV-inverter connected to phase $\phi \in \cP_i$ of node $i \in \cN_D$. Lines $(m,n) \in \cE$ are still modeled as $\pi$-equivalent components~\cite[Ch. 6]{Kerstingbook} and the $|\cP_{mn}| \times |\cP_{mn}|$ phase impedance and shunt admittance matrices are denoted as $\bZ_{mn} \in \mathbb{C}^{|\cP_{mn}| \times |\cP_{mn}|}$ and $\bY_{mn}^{(s)} \in \mathbb{C}^{|\cP_{mn}| \times |\cP_{mn}|}$, respectively.  Three- or single-phase transformers (if any) are modeled as series components with transmission parameters that depend on the connection type~\cite[Ch.~8]{Kerstingbook},~\cite{Paudyal11}. 

A prototypical non-convex OPF formulation can be readily obtained by enforcing the balance constraints and the voltage regulation constraints on a per-node and per-phase basis~\cite{Paudyal11,Dallanese-TSG13,Gan-PSCC14}; and, by properly augmenting the cost function to account for (cost of) power losses and power injections over all phases. To develop an SDP relaxation of the non-convex three-phase OPF problem,  consider re-defining the vector of voltages $\bv$ as $\bv := [\bv_0^\sfT,\bv_1^\sfT, \ldots,\bv_{N}^\sfT]^\sfT$, where $\bv_i := [\{V_i^{\phi}\}_{\phi \in \cP_i}]^\sfT$ is a $|\cP_i| \times 1$ vector collecting the voltages on the phases of node $i \in \cN$. Similarly, vector $\bi$ now collects the currents injected in all nodes and phases; that is, $\bi := [\bi_0^\sfT,\bi_1^\sfT, \ldots,\bi_{N}^\sfT]^\sfT$. As shown in Section~\ref{sec:ProblemFormulation}, Ohm's law and Kirchhoff's current law can be captured by the linear equation $\bi = \bY \bv$ where, in this case, the network admittance matrix $\mathbf{Y}$ has dimensions $\sum_{i\in\cN} |\cP_i| \times \sum_{i\in\cN} |\cP_i|$, and its entries are computed based on the system topology and the lines matrices $\{\bZ_{ij}, \bY_{ij}^{(s)}\}_{(i,j) \in \cE}$ as specified in~\cite{Dallanese-TSG13,Gan-PSCC14} and~\cite{Robbins15}. 
To express voltage magnitudes and powers as linear functions of the outer-product matrix $\bV := \bv \bv^\sfH$, define matrix $\bY_i^{\phi} := \bar{\be}_i^{\phi} (\bar{\be}_i^{\phi})^{\sfT} \bY$  per node $i$ and phase $\phi$, where $\bar{\be}_i^{\phi} := [\mathbf{0}_{|\cP_0|}^{\sfT},\ldots,\mathbf{0}_{|\cP_{i-1}|}^{\sfT},\be_{\cP_i}^{\phi, \sfT},\mathbf{0}_{|\cP_{i+1}|}^{\sfT},\ldots,\mathbf{0}_{|\cP_N|}^{\sfT}]^{\sfT}$, and $\{\be_{\cP_i}^\phi\}_{\phi \in \cP_n}$ denotes the canonical basis of $\mathbb{R}^{|\cP_i|}$.  Next, per node $i \in \cN$ and phase $\phi \in \cP_i$, define the Hermitian matrices $\bPhi_{i}^{\phi} := \frac{1}{2} (\bY_i^{\phi} + (\bY_i^{\phi})^\sfH)$,  $\bPsi_{i}^{\phi} := \frac{j}{2} (\bY_i^{\phi} - (\bY_i^{\phi})^\sfH)$, and $\bUpsilon_{i}^{\phi} := \bar{\be}_i^{\phi} (\bar{\be}_i^{\phi})^{\sfT}$. Then, the net real and reactive powers injected at node $i$ and phase $\phi$ can be expressed as $\trace(\bPhi_{i}^{\phi} \bV) = P_{i}^{\phi} - P_{\ell,i}^{\phi}$ and $\trace(\bPsi_{i}^{\phi} \bV) = Q_{i}^{\phi} - Q_{\ell,i}^{\phi}$, respectively, whereas the squared voltage magnitude at the same node and phase reads $\trace(\bUpsilon_i^{\phi} \bV)  = |V_{i}^{\phi}|^2$.

Using these definitions, the OPF problem can be formulated as:
\begin{subequations}
\label{eq:p1-threephase}
\begin{align}
  \min_{\bV \succeq \mathbf{0}, \{\bu_i^\phi \in \cY_i^\phi \}}& \,\, H(\bV) + \sum_{i \in \cN_D} \sum_{\phi \in \cP_i} G_i^{\phi}(\bu_i^\phi) \label{eq:p1-threephase-cost} \\
& \hspace{-1.8cm} \mathrm{subject~to}  \,\, \bh_i^\phi(\bV)  =  \bu_i^\phi - \bd_i^\phi  , \hspace{.45cm} \forall \, i \in \cN_D , \phi \in \cP_i\label{eq:p1-threephase-eq} \\
& \hspace{-.1cm} \bh_i^\phi(\bV)  = - \bd_i^\phi  ,  \hspace{1.0cm} \forall \, i \in \cN_O , \phi \in \cP_i \label{eq:p1-threephase-eq-load} \\
& \hspace{-.9cm} V_{\mathrm{min}}^2 \leq \trace(\bUpsilon_i^\phi \bV)  \leq V_{\mathrm{max}}^2,  \hspace{.2cm} \forall \, i \in \cN , \phi \in \cP_i \label{eq:p1-threephase-eq-load} \\
&  \mathrm{rank}(\bV) = 1 \label{eq:p1-threephase-rank} 
\end{align}
\end{subequations}
where $ \bh_i^\phi(\bV)  = [\trace(\bPhi_i^\phi \bV), \trace(\bPsi_i^\phi \bV)]^\sfT$, $\bu_i^\phi = [P_i^\phi, Q_i^\phi]^\sfT$, and $\bd_i^\phi = [P_{\ell,i}^\phi, Q_{\ell,i}^\phi]^\sfT$ [cf.~\eqref{eq:correspondence_h}]. An SDP relaxation of problem~\eqref{eq:p1-threephase} can be obtained by discarding the rank constraint~\eqref{eq:p1-threephase-rank}.

The procedure outlined in Section~\ref{sec:controllerSynthesis} can  be utilized  to 
synthesize controllers for the PV-inverters that solve~\eqref{eq:p1-threephase} in a recursive manner. To this end, if suffices to dualize balance constraints~\eqref{eq:p1-threephase-eq} to form the (partial) Lagrangian~\eqref{eq:lagrangian}, and follow steps~\eqref{eq:dualsubgradient_system_as}. In particular, the resultant distributed algorithm entails the following operations:
\begin{subequations}
\label{eq:dualsubgradient_system_as_phase}
\begin{align}
&  \hspace{.15cm} \by_i^\phi[t_{k}] =  \br_i^\phi  \Big( \bx_i^\phi (t_{k}), \bd_i^\phi  \Big)  \label{eq:sys1-obs_sys_as_phase}  \\
& \hspace{-.2cm}  \blambda_i^\phi [t_{k+1}]   = \blambda_i^\phi [t_{k}]  + \alpha_{k+1} \Big( \bh_i^\phi (\bV[t_{c(k)}]) - \by_i^\phi [t_{k}] + \bd_i^\phi   \Big)   \label{eq:dual_ascent_sys_as_phase} \\
& \hspace{-.2cm}  \bu_i^\phi  [t_{k+1}]  = \arg \min_{\bu_i^\phi  \in \cY_i^\phi } G_i^\phi (\bu_i^\phi ) - (\blambda_i^\phi[t_{k+1}])^\sfT \bu_i^\phi  \, \label{eq:primal_y_sys_as_phase} 
\end{align}
\begin{align}
\hspace{-.2cm} \bV[t_{c(k)}]  = \arg \min_{\bV \in \bar{\cV}} \, H(\bV) +  \sum_{i \in \cN_D} \sum_{\phi \in \cP_i} (\blambda_i^\phi[t_{c(k)}])^\sfT \, \bh_i^\phi (\bV) \label{eq:primal_V_sys_as_phase}  \hspace{-.2cm} 
\end{align}
\end{subequations}
where~\eqref{eq:sys1-obs_sys_as_phase}--\eqref{eq:primal_y_sys_as_phase} are performed at each PV-inverter connected to phase $\phi$ of node $i$, and~\eqref{eq:primal_V_sys_as_phase} is carried out at the DSO.  In~\eqref{eq:primal_V_sys_as_phase}, the set $\bar{\cV}$ is defined as   
\begin{align}
& \bar{\cV} :=  \{\bV: \bV \succeq \mathbf{0}, V_{\mathrm{min}}^2 \leq \trace(\bUpsilon_i^\phi \bV)  \leq V_{\mathrm{max}}^2 \,\forall \, i \in \cN  \, \nonumber  \\
& \mathrm{~and~}  \bh_i^\phi(\bV)  = - \bd_i^\phi, \,  \forall \phi \in \cP_i, i \in \cN_O \} \, .  \label{eq:V_region_phase}  
\end{align} 
It can be readily shown that the convergence claims of \emph{Theorem~\ref{thm:convergence_as}} carry over to the multi-phase setup.
 
\bibliographystyle{IEEEtran}
\bibliography{biblio.bib}

\end{document}